\documentclass[letterpaper,11pt,leqno,oneside]{article}
\usepackage{color}
\usepackage{bm}
\usepackage{exscale}
\usepackage{amsmath}
\usepackage{amsfonts}
\usepackage{stmaryrd}
\usepackage{amscd}
\usepackage{graphicx}
\usepackage{amsxtra}
\usepackage{amssymb}
\usepackage{theorem}
\usepackage[final]{epsfig}
\usepackage{eqnarray}
\usepackage{hyperref}
\newtheorem{proposition}{Proposition}[subsection]
\newtheorem{definition}[proposition]{Definition}

\newtheorem{lemma}[proposition]{Lemma}

{\theorembodyfont{\rmfamily}\newtheorem{remark}[proposition]{Remark}}
\newtheorem{theorem}[proposition]{Theorem}
\newtheorem{corollary}[proposition]{Corollary}

{\theorembodyfont{\rmfamily}\newtheorem{example}[proposition]{Example}}

\newfont{\abc}{cmtt10 scaled 1200}

\def\R{\mathbb{R}}

\def\Z{\mathbb{Z}}

\def\U{\mathbb{U}}
\def\D{\mathbb{D}}

\def\E{\mathbb{E}}
\def\U{\mathbb{U}}

\def\I{\mathbb{I}}
\def\ve{\varepsilon}

\def\ra{\rightarrow}
\def\cs{\symbol{35}}
\def\p{\partial}
\def\qed{\hfill $\Box$ \\}

\def\mm{\mbox}
\def\v{= \emptyset}
\def\n{\neq \emptyset}
\def\D{\mathbf{ID}}
\def\M{\mathbb{A}}

\def\sks{$|A|$-skins}

\def\bp{\langle A \rangle}
\def\db{d^{\,\flat}}

\def\si{$\mathcal{S}$}

\def\llb{\llbracket}
\def\rrb{\rrbracket}

\def\sks{$|A|$-skins}
\def\bp{\langle A \rangle}

\begin{document}
\vspace*{0cm}

\begin{center}\Large{\bf{Hyperbolic Unfoldings of Minimal Hypersurfaces}}\\
\bigskip
\large{\bf{Joachim Lohkamp}}
\end{center}

\noindent Mathematisches Institut, Universit\"at M\"unster, Einsteinstrasse 62, 48149 M\"unster, Germany\\
{\small{\emph{e-mail: j.lohkamp@uni-muenster.de}}}\\

\noindent\textbf{Abstract:} We study the intrinsic geometry of area minimizing hypersurfaces from a new point of view by relating this subject to quasiconformal geometry. Namely, for any such hypersurface $H$ we define and construct a so-called \si-\emph{structure}. This new and natural concept reveals some unexpected geometric and analytic properties of $H$ and its singularity set $\Sigma$. Moreover, it can be used to prove the existence of hyperbolic unfoldings of $H\setminus \Sigma$. These are canonical conformal deformations of $H \setminus \Sigma$ into complete Gromov hyperbolic spaces of bounded geometry with Gromov boundary homeomorphic to $\Sigma$. These new concepts and results naturally extend to the larger class of almost minimizers.\\



\tableofcontents

%
%
%
%
%%%%%%%%%%%%%%%%%%%%%%%%%%%%%%%%%
%%%%%%%%%%%%%%%%%%%%%%%%%%%%%%%%%
%%%%%%%%%%%%%%%%%%%%%%%%%%%%%%%%%
\setcounter{section}{1}
\renewcommand{\thesubsection}{\thesection}
\subsection{Introduction} \label{int}
%%%%%%%%%%%%%%%%%%%%%%%%%%%%%%%%%
%%%%%%%%%%%%%%%%%%%%%%%%%%%%%%%%%
%%%%%%%%%%%%%%%%%%%%%%%%%%%%%%%%%

%
\noindent Let $M^{n+1}$ be a smooth compact manifold, and $H^n\subset M^{n+1}$ be an area minimizing hypersurface with singularity set $\Sigma\subset H$. It is known that $\Sigma$ is a potentially complicated compact set of Hausdorff-dimension $\le n-7$ with some serious impact also on $H \setminus \Sigma$. The second fundamental form $A_H$ and its norm $|A_H|$ diverge towards $\Sigma$. The open manifold $H \setminus \Sigma$ collapses while we approach $\Sigma$ so that even the topology of arbitrarily small balls in $H$, around a given singular point, can be highly non-trivial.\\

To manage this complex situation we establish structures on $H \setminus \Sigma$ which help to understand the geometric analysis of and also on $H \setminus \Sigma$ without using the structure
of $\Sigma$.\\

\noindent\textbf{\si-structures.}\, The key idea of this paper is to introduce on these hypersurfaces natural distance and size concepts, the \emph{\si-structures}, which measure also the curvature of $H$. For instance, we get the \si-distance $\delta_{\bp}$ which measures a generalized form of distance to the singular set and which commutes with blow-ups around singular points. For the ordinary metric distance this commutativity fails. The option to employ blow-ups is one of the reasons why \si-structures simplify the study of geometric analysis on $H \setminus \Sigma$ near $\Sigma$.\\

\noindent$\bullet$\, Regarding $\Sigma$ as the boundary of the open manifold $H \setminus \Sigma$, \si-structures unravel some global boundary regularity for $H\setminus\Sigma$, namely its \emph{uniformity} and the even stronger \emph{\si-uniformity}. The uniformity concept arose from the study of Euclidean domains with highly irregular boundary, but which still retain many geometro-analytic properties of smooth domains.\\

\noindent$\bullet$\, \si-uniformity also takes the curvature degeneration of $H \setminus \Sigma$ towards $\Sigma$ into account (though it remains a non-trivial concept even when $\Sigma \v$). This, and not merely uniformity, is the essential tool to prove existence of \emph{hyperbolic unfoldings} of $H \setminus \Sigma$. These are canonical conformal deformations of $H \setminus \Sigma$ into complete \emph{Gromov hyperbolic spaces} of \emph{bounded geometry}. Moreover, the \emph{Gromov boundary} $\p_G(H \setminus \Sigma)$ of such an unfolding is just the singular set, i.e., it is homeomorphic to $\Sigma \subset H$.\\

Basic ingredients to derive the \si-uniformity and the existence of hyperbolic unfoldings are the isoperimetric inequality and the regularity theory for area minimizers. A further distinctive property we use only holds in the case of  \emph{hypersurfaces}. Namely, their tangent cones at singular points are also embedded singular hypersurfaces.\\

 \noindent$\bullet$\,  Our results equally apply to the larger class of \emph{almost} minimizers. They can the characterized as possibly singular hypersurfaces which asymptotically look like area minimizers when we approach their singular set. This class includes hypersurfaces with prescribed mean curvature or obstacles, and also cases not arising from variational problems like hypersurfaces evolving under geometric flows or occurring as horizons of black holes in general relativity. In a similar vein, we can treat (almost) area minimizers with boundaries (solving a \emph{Plateau problem}). However, to keep the arguments easier to follow we confine ourselves to the more familiar case of area minimizers in the main text and postpone their extension to almost minimizers to Appendix A.III.\\

\noindent\textbf{Typical applications.}\,
The combination of hyperbolicity and bounded geometry simplifies dramatically the geometric analysis on the hyperbolic unfolding of $H \setminus \Sigma$. For instance, building on Ancona's work~\cite{An1},~\cite{An2} we can start to work out the potential theory of many naturally defined elliptic operators. These results can then be referred back to the original space $H \setminus \Sigma$. In by-passing the difficult internal structure of $\Sigma$, hyperbolic unfoldings become a versatile tool for the very delicate geometric analysis on singular area minimizers.\\

%
%
%
%%%%%%%%%%%%%%%%%%%%%%%%%%%%%%%%%
\subsubsection{Basic Notations}\label{notation}
%%%%%%%%%%%%%%%%%%%%%%%%%%%%%%%%%
In this paper $H^n$ denotes a connected integer multiplicity rectifiable current of dimension $n \ge 2$ which sits inside some complete, smooth Riemannian manifold $(M^{n+1}, g_M)$. We briefly refer to such a current $H$ as an \emph{area minimizer} when it is a locally mass minimizing. By $\Sigma_H$, or simply $\Sigma$ if there is no risk of confusion, we denote the set of singular points of $H$. (For the convenience of the reader we recall some facts from geometric measure theory in Appendix A.) For a minimal cone $C$ we write the singular set $\sigma_C$ as a hint that we think of them as tangential spaces. The upper/lower case notation inspired from the case Lie groups and their Lie algebras.\\

For any $A \subset X$, in a metric space $(X,d_X)$, the distance to $A$ is denoted by $dist_{d_X}(\cdot, A)$. By a curve we mean a continuous map $\gamma:[a,b] \ra X$, $a<b$. Its length $l_{d_X}(\gamma)$ is defined by {\small $l_{d_X}(\gamma):= \sup\{\sum_{i=0,..,N}d_X(\gamma(t_{i-1}),\gamma(t_i))\,\Big|\, \mm{partitions }a=t_0 \le t_1 \le ... \le t_N = b\}$}.  $\gamma$ is \emph{rectifiable}  if $l_{d_X}(\gamma) < \infty$. $X$ is  \emph{rectifiably connected} if any two $p,q \in X$ can be joined by a rectifiable curve. For a Riemannian manifold $(X,g_X)$  we also directly use $g_X$ as an index  in place of its associated metric $d_X$. When there are no ambiguities we usually omit these indices.\\

 The Riemannian metric on $H$ induced for its embedding $H \subset M$ is denoted by $g_H$. For $(H,g_H)$ viewed as a metric space we refer to the induced distance function $d_{g_H}(p,q)$ for $p,q \in H$ as the  \emph{intrinsic distance}, whereas $d_{g_M}(p,q)$  is the \emph{extrinsic distance} relative $M$.
  For $\lambda>0$ we let $\lambda\cdot M$ denote the conformally rescaled Riemannian manifolds $(M,\lambda^2 \cdot  g)$. In the sequel, we shall consider the following classes of \emph{complete} area minimizers:
\begin{description}
  \item[${\cal{H}}^c_n$:] $H^n \subset M^{n+1}$ is a compact embedded hypersurface without boundary.
  \item[${\cal{H}}^{\R}_n$:] $H^n \subset\R^{n+1}$ is a complete hypersurface in flat Euclidean space $(\R^{n+1},g_{\R^{n+1}})$ with $0\in H$ and which is an oriented boundary of some open set in $\R^{n+1}$.
  \item[${\cal{H}}_n$:] ${\cal{H}}_n:= {\cal{H}}^c_n \cup {\cal{H}}^{\R}_n$ and we set ${\cal{H}} :=\bigcup_{n \ge 1} {\cal{H}}_n$.
\end{description}
\begin{remark} \,
Any current in ${\cal{H}}_n$ can be locally decomposed into (locally disjoint) oriented minimal boundaries of open sets, cf.\ Appendix A, Propositions~\ref{dic} and \ref{dicc} as well as~\cite[4.5.17]{F1}, \cite[Chapter 37]{Si1} and \cite{Si2}. Consequently, we may assume that $H$ is \emph{locally} an oriented boundary of an open set in $M$.
\end{remark}
We shall consider the following larger classes of \emph{almost} minimizers cf. Appendix A. II:
\begin{description}
  \item[${\cal{G}}^c_n$:] $H^n \subset M^{n+1}$ is a compact connected  almost minimizer. We set ${\cal{G}}^c :=\bigcup_{n \ge 1} {\cal{G}}^c_n.$
  \item [${\cal{G}}_n$:] ${\cal{G}}_n := {\cal{G}}^c_n \cup {\cal{H}}^{\R}_n$ and ${\cal{G}} :=\bigcup_{n \ge 1} {\cal{G}}_n$. We notice ${\cal{H}}^c_n \varsubsetneq {\cal{G}}^c_n$ and ${\cal{H}}_n \varsubsetneq {\cal{G}}_n$.
\end{description}
\begin{remark} \,
Even if one is merely interested in ${\cal{H}}^c_n$ or ${\cal{G}}^c_n$ it is important to include ${\cal{H}}^\R_n$ in all arguments. This way we get spaces ${\cal{H}}_n$ and ${\cal{G}}_n$ which are closed under \emph{blow-ups} of area minimizers, cf.\ Example~\ref{ex1}, and we can use compactness results, on the space ${\cal{H}}^{\R}_n$, in the study of $H \in {\cal{H}}^c_n$ or ${\cal{G}}^c_n$ near $\Sigma_H$.
\end{remark}
%
%
%
%%%%%%%%%%%%%%%%%%%%%%%%%%%%%%%%%
\subsubsection{Main definitions and results}\label{sor}
%%%%%%%%%%%%%%%%%%%%%%%%%%%%%%%%%
There are two lines of results. The first line concerns the global boundary regularity of the open manifold $H\setminus\Sigma$, the so-called \emph{\si-uniformity}. The second line explores \emph{hyperbolic unfoldings} and reveals the hyperbolic nature of area minimizers.\\

\noindent\textbf{\si-transforms.} \,
The key tool for proving these results is the construction of an \si-transform $\bp_H$ on an area minimizer $H \in {\cal{H}}$ (or more generally on some $H \in {\cal{G}}$). It results from a particular way of merging the induced Riemannian metric $g_H$ and the second fundamental form $A=A_H$ into a scalar function $\bp_H$. The characteristic property is that its level sets $\bp^{-1}_H(c)$ can be thought of as regularizing  ``membranes'' or ``skins'' spanned over the barely controlled level sets $|A|^{-1}(c)$. (This concept arose from the idea to describe well-controlled and naturally defined domains with closure in  $H \setminus \Sigma$ so that their (elliptic) analysis efficiently approximates the global analysis on $H \setminus \Sigma$.) The label \si\ then stands for both \emph{skin systems} and the resulting \emph{strong or super}-uniformity, namely \si-uniformity.\\

We will construct a concrete \si-transform below. However, different \si-transforms still share some basic properties/axioms. Our applications and arguments only employ these few properties of  a concrete \si-transform so that we give the following axiomatic definition.

\begin{definition}[\si-transforms]\label{def1} We call an assignment $\bp$ that associates with any $H \in {\cal{G}}$ a function $\bp_H:H \setminus \Sigma_H\to\R$ an \textbf{\si-transform} provided it satisfies the following axioms:
\begin{description}
  \item[(S1)] \emph{\textbf{Trivial Gauge}} \,
  If $H \subset M$ is totally geodesic, then $\bp_H \equiv 0$.
  \item[(S2)] \emph{\textbf{\si-Properties}} \,
  If $H$ is not totally geodesic, then the level sets $\M_c:= \bp_H^{-1}(c)$, for $c>0$, we call the $|A|$-\textbf{skins}, surround the level sets of $|A|$:
  \[
  \bp_H>0, \bp_H \ge |A_H| \mm{ and } \bp_H(x) \ra \infty, \mm{ for } x \ra p \in \Sigma_H.
  \]
Like $|A_H|$, $\bp_H$ anticommutes with scalings, i.e., $\bp_{\lambda \cdot H} \equiv \lambda^{-1} \cdot  \bp_{H}$ for any $\lambda >0$.
  \item[(S3)] \emph{\textbf{Lipschitz regularity}} \,
  If $H$ is not totally geodesic, and thus $\bp_H>0$, we define the
  \[
  \mm{ \textbf{\si-distance} } \delta_{\bp_H}:=1/\bp_H.
  \]
  This function is
  $L_{\bp}$-Lipschitz regular for some constant $L_{\bp}=L(\bp,n)>0$, i.e.,
  \[
  |\delta_{\bp_H}(p)- \delta_{\bp_H}(q)|   \le L_{\bp} \cdot d_{g_H}(p,q) \mm{ for any } p,q \in  H \setminus \Sigma \mm{ and any } H \in {\cal{G}}_n.
  \]
If $H$ is totally geodesic, and thus $\bp_H \equiv 0$, we set $\delta_{\bp_H}\equiv\infty$ and use the convention $|\delta_{\bp_H}(p)- \delta_{\bp_H}(q)|\equiv0$.
  \item[(S4)]  \emph{\textbf{Naturality}} \,
  If $H_i \in {\cal{H}}_n$, $i \ge 1$, is a sequence converging* to the limit space $H_\infty \in {\cal{H}}_n$, then $\bp_{H_i}\overset{C^\alpha}  \longrightarrow {\bp_{H_\infty}}$ for any $\alpha \in (0,1)$. For general $H \in {\cal{G}}_n$, this holds for blow-ups:  $\bp_{\tau_i \cdot H} \overset{C^\alpha}  \longrightarrow {\bp_{H_\infty}}$, for any sequence $\tau_i \ra \infty$ so that  $\tau_i \cdot H \ra H_\infty \in {\cal{H}}^\R_n$.
\end{description}
\end{definition}
*For the precise notions of convergence we use here, see Section~\ref{nat} and Appendix A.II, III. To simplify notation, we omit the index $H$ in $\bp_H$ and $\delta_{\bp_H}$ if there is no risk of confusion.

\begin{remark} \,1. If $H \in {\cal{G}}$ and $H \subset M$ is totally geodesic, it is not hard to see that $\Sigma \v$, cf.\ Corollary~\ref{t3} from Appendix A. In this paper, the totally geodesic hypersurfaces are the \emph{trivial} cases: many results either obviously hold or they degenerate to conventions.\\
2. The only Lipschitz regular $\delta_{\bp}$ can be approximated by some Whitney type $C^\infty$-smoothing $\delta_{\bp^*}$ satisfying (S1)-(S3) with $c_1 \cdot \delta_{\bp}(x) \le \delta_{\bp^*}(x)  \le c_2 \cdot \delta_{\bp}(x)$, for some constant  $c_1>0$, cf. Appendix B, Proposition~\ref{smsk}.
\end{remark}

To prove the mere existence of \si-structures we use an interpolation between the functions $|A|$ and $1/dist_{g_H}(x,\Sigma)$. These so-called \emph{metric \si-structures} which result from this procedure have some additional properties and account for our basic intuition on \si-transforms.

\begin{theorem}[Metric \si-transforms]\label{thm1} \,
There is a family of \si-transforms, $\bp_\alpha$, $\alpha >0$, we call the \textbf{metric \si-transforms}, with the following properties:
 \begin{itemize}
  \item \, The \sks\ $\M_c$ of $\bp_\alpha$ bound the outer $\alpha/c$-distance collar of $|A|^{-1}[c,\infty)$ in $H$.
    \item \,   $\bp_{\alpha}(x) \ra |A|(x)\mm{ in } L_{loc}^\infty,\mm{ for }\alpha \ra 0, \mm{ on } H \setminus \Sigma$
        \item \,   $1/\alpha \cdot \bp_{\alpha}(x) \ra 1/ dist_{g_H}(x,\Sigma)\mm{ in } L_{loc}^\infty,\mm{ for }\alpha \ra \infty, \mm{ on } H \setminus \Sigma$.
  \end{itemize}
\end{theorem}

We will prove Theorem~\ref{thm1} in Section~\ref{st}.

\begin{remark} \,
The limit cases $|A|$ and $1/ dist_{g_H}(\cdot,\Sigma)$ are no longer \si-transforms. In general $|A|^{-1}(c)\cap \Sigma \n$ and there is no uniform Lipschitz bound for $|A|^{-1}$. Thus $|A|$ violates (S2) and (S3). On the other hand, there is no $c>0$ so that $c/ dist_{g_H}(\cdot,\Sigma) \ge |A|$. Furthermore, there is no proper correlation between the singularities of converging sequences in ${\cal{H}}$ and of their limit (consider for instance a family of smooth area minimizers converging to a singular one). Thus $1/ dist_{g_H}(\cdot,\Sigma)$ violates (S2) and (S4).
\end{remark}

For the remainder of this introduction (and all later applications) we consider a general \si-transform $\bp$ satisfying the axioms of Def.~\ref{def1}. The choice of a different \si-transform merely changes the global constants in the statements below.\\

\noindent\textbf{\si-uniformity.} \,
For the analysis and geometry near the boundary $\Sigma$ it is crucial to quantify the approachability of $\Sigma$ from within $H \setminus \Sigma$. For our purposes we need a global boundary regularity condition which ensures non-tangential accessibility of $\Sigma$. An appropriate starting point is the notion of uniform space, see for instance~\cite{BHK}, \cite{He}.

\begin{definition}[Uniform spaces]\label{ud} \, Let $(X,d)$ be a non-complete, locally compact and complete, rectifiably connected metric space. We denote its metric completion by $\overline{X}$ and define its boundary by $\p X:= \overline{X} \setminus X$. $(X,d)$ is a \textbf{c-uniform space}, or a \textbf{uniform space} for short, if there is a $c \in \R^{\ge 1}$, so that any two points $p,q \in X$ can be joined by a
\textbf{c-uniform curve}. That is a rectifiable path $\gamma: [a,b] \ra X$, for some $a <b$, from $p$ to $q$ so that
\begin{itemize}
  \item \emph{\textbf{Quasi-geodesic:}} \, $l(\gamma_{p,q})\le c\cdot d(p,q)$.
  \item \emph{\textbf{Twisted double cones:}} \, Let $l_{min}(\gamma_{p,q}(z)):=$ minimum of the lengths of the two subcurves of $\gamma_{p,q}$ from $p$ to $z$ and from $q$ to $z$. Then
  \[
  l_{min}(\gamma_{p,q}(z)) \le c \cdot dist(z,\p X), \mm{ for any } z \in \gamma_{p,q}.
  \]
\end{itemize}
\end{definition}

We demonstrate a stronger form of this uniformity of $H \setminus \Sigma$, its \si-uniformity, for any $H \in {\cal{H}}_n$. The new concept also naturally extends to the regular case where $\Sigma \v$. In Sections \ref{rn} and \ref{inta} we will prove the

\begin{theorem}[\si-Uniformity of $\mathbf{H \setminus \Sigma}$]\label{thm2} \,
Let $H \in {\cal{G}}$ be a hypersurface with (possibly empty) singular set $\Sigma=\Sigma_H$.
\begin{enumerate}
  \item$H \setminus \Sigma$ and $H$ are \textbf{rectifiably connected}. In particular, any compact $H \in {\cal{G}}^c$ has a finite intrinsic diameter: $diam_{g_H}H <\infty.$
  \item There exists $c>0$ such that $H \setminus \Sigma$ is a \textbf{c-\si-uniform space}, or \textbf{\si-uniform space} for short. This means that any pair $p,q \in H \setminus \Sigma$ can be joined by a \textbf{c-\si-uniform curve} in $H \setminus \Sigma$, i.e., a rectifiable curve $\gamma_{p,q}: [a,b] \ra H \setminus \Sigma$ with $\gamma_{p,q}(a)=p$, $\gamma_{p,q}(b)=q$ and such that the following properties hold.
  \begin{itemize}
    \item \emph{\textbf{Quasi-geodesic:}} \, $l_{g_H}(\gamma)  \le c \cdot  d_{g_H}(p,q).$
    \item \emph{\textbf{Twisted double \si-cones:}} \, $l_{min}(\gamma_{p,q}(z)) \le c \cdot \delta_{\bp}(z)$ for any $z \in \gamma_{p,q}$.
  \end{itemize}
\end{enumerate}
\end{theorem}

\begin{remark} \,
As a first application we see that \si-uniformity implies uniformity of $H \setminus \Sigma$ if $\Sigma \n$, a result which would be hard to derive directly. Indeed, the Lipschitz condition (S3) implies $\delta_{\bp}(x) \le L  \cdot dist_{g_H}(x,\Sigma)$ for any $x \in H \setminus \Sigma$. For totally geodesic $H$ the result holds trivially since $\delta_{\bp} \equiv +\infty$ and either $H$ is compact and smooth or a Euclidean hyperplane.
\end{remark}

\begin{remark} \,
For $H\in{\cal{H}}^{\R}_n$ the \si-uniformity parameter $c$ depends only on $n$.
\end{remark}

\noindent\textbf{Hyperbolic geometry on $H \setminus \Sigma$.} \,
By results of Gehring and Osgood~\cite{GO} and Bonk, Heinonen and Koskela~\cite{BHK} uniform spaces are Gromov hyperbolic. More concretely, for any \emph{singular} hypersurface $H \in {\cal{G}}$ we can define the \emph{\textbf{quasi-hyperbolic metric}}
\[
k_{H \setminus \Sigma}(x,y) := \inf \Bigl \{\int_\gamma 1/dist_{g_H}(\cdot,\Sigma_H)  \, \, \Big| \, \gamma  \subset H \setminus \Sigma \mbox{ rectifiable curve joining }  x \mbox{ and } y  \Bigr \},
\]
for $x$, $y\in H \setminus \Sigma$. Moreover, $(H \setminus \Sigma, k_{H \setminus \Sigma})$ is a \emph{complete, Gromov hyperbolic and visual metric space} (see Section~\ref{discrete} for a definition of these concepts).\\

However, uniformity does not take into account the smooth but highly curved regions of $H$ near $\Sigma$, that is, the geometry of $(H \setminus \Sigma, k_{H \setminus \Sigma})$ might \emph{not} be bounded. This, on the other hand, makes the analysis of elliptic operators with respect to $k_{H \setminus \Sigma}$ quite a subtle endeavor. Moreover, like the singular set, $k_{H \setminus \Sigma}$ may change drastically even after small deformations of $H$. This limits the use of this metric in blow-up and compactness arguments.\\

At any rate, from the viewpoint of the discerned \si-uniformity of $H \setminus \Sigma$ there is a more versatile and natural hyperbolic geometry on $H \setminus \Sigma$, the \textbf{\si-metric} $d_{\bp_H}$, which resolves the issues with $k_{H \setminus \Sigma}$. This metric is defined by
\[
d_{\bp}(x,y) := \inf \Bigl  \{\int_\gamma  \bp \, \, \Big| \, \gamma   \subset  H \setminus \Sigma\mbox{ rectifiable curve joining }  x \mbox{ and } y  \Bigr \}
\]
for $x$, $y\in H \setminus \Sigma$. Recall that $\bp_H=1/\delta_{\bp_H}$ to see the analogy with the definition of $k_{H \setminus \Sigma}$.  Here and in what follows we drop the index $H$ when it is known from the context. In terms of our general conventions we have $d_{\bp} \equiv d_{\bp^2 \cdot g_H}$.
In addition, the \si-metric is also defined for \emph{smooth} $H$ where $\Sigma\v$.
The \si-uniformity is the key ingredient in proving the following result for \si-metrics in Section~\ref{gh}. It is our main hyperbolization theorem.

\begin{theorem}[Conformal hyperbolic unfoldings]\label{thm3} \,
For any non-totally geodesic hypersurface $H \in {\cal{G}}$, the \si-metric $d_{\bp}$ has the following properties:
\begin{itemize}
  \item The metric space $(H \setminus \Sigma, d_{\bp})$ and its \textbf{quasi-isometric Whitney smoothing}, i.e.\ the smooth Riemannian manifold $(H \setminus \Sigma, d_{\bp^*}) = (H \setminus \Sigma, 1/\delta_{\bp^*}^2 \cdot g_H)$, are \textbf{complete Gromov hyperbolic spaces} with \textbf{bounded geometry}.
  \item $d_{\bp}$ is natural. That is, the assignment of $d_{\bp_H}$ to $H$ commutes with the compact convergence of the regular portions of the underlying area minimizers.
\end{itemize}
The spaces $(H \setminus \Sigma, d_{\bp})$ and $(H \setminus \Sigma, d_{\bp^*})$ are conformally equivalent to the original space $(H \setminus \Sigma, g_H)$. We refer to both these spaces as \textbf{hyperbolic unfoldings} of $(H \setminus \Sigma, g_H)$.
\end{theorem}

\begin{remark} \,
Again, for $H\in{\cal{H}}^{\R}_n$ the hyperbolicity and the bounded geometry parameters depend only on $n$. For totally geodesic $H \in {\cal{G}}$, both $(H \setminus \Sigma, d_{\bp})$ and $(H \setminus \Sigma, d_{\bp^*})$ are still well-defined but they are one-point spaces since $\bp \equiv 0$ (and $\Sigma \v$ cf. ~\ref{t3}.) In all other cases, $\bp>0$ and, hence,  $(H \setminus \Sigma, d_{\bp})$ and $(H \setminus \Sigma, d_{\bp^*})$ are homeomorphic  to $(H \setminus \Sigma, g_H)$.
\end{remark}

Finally, we consider the Gromov boundaries of these hyperbolic spaces. Let us denote the one-point compactification of a hypersurface $H \in {\cal{H}}^{\R}_n$ by $\widehat{H}$. For the singular set $\Sigma_H$ of some $H \in {\cal{H}}^{\R}_n$ we \emph{always} add $\infty_H$  to $\Sigma$ and define $\widehat\Sigma:=\Sigma \cup \infty_H$ (note that $\Sigma$ could already be compact). For  $H \in {\cal{G}}^{c}_n$ we set $\widehat H=H$ and $\widehat\Sigma=\Sigma$. In Section \ref{grch} we prove the following theorem rendering $\widehat\Sigma$ as the Gromov boundary of the hyperbolic unfoldings of $(H \setminus \Sigma, g_H)$.

\begin{theorem}[Gromov boundary of $H \setminus \Sigma$]\label{thm4} \,
For any non-totally geodesic  $H \in {\cal{G}}$ the identity map on $H \setminus \Sigma$ extends to homeomorphisms between the one-point compactification $\widehat{H}$ and the Gromov compactifications $\overline{X}_G$ of $X=(H \setminus \Sigma,d_{\bp})$, $(H \setminus \Sigma,d_{\bp^*})$ and $(H \setminus \Sigma, k_{H \setminus \Sigma})$:
\[
\widehat{H}\cong\overline{(H \setminus \Sigma,d_{\bp})}_G \cong \overline{(H \setminus \Sigma,d_{\bp^*})}_G \cong \overline{(H \setminus \Sigma, k_{H \setminus \Sigma})}_G,
\]
where $ \cong$ means homeomorphic. In particular, we find for the Gromov boundaries $\p_G(X)$:
\[
\widehat{\Sigma} \cong\p_G(H \setminus \Sigma,d_{\bp}) \cong \p_G(H \setminus \Sigma,d_{\bp^*}) \cong \p_G(H \setminus \Sigma, k_{H \setminus \Sigma}).
\]
\end{theorem}
\begin{remark} \,
For \emph{smooth} hypersurfaces i.e., for $\Sigma \v$, this reads as follows. In the case where $H\in{\cal{G}}^c_n$ the Gromov boundary is empty and  the hyperbolic unfolding is again a compact manifold without boundary. For $H\in{\cal{H}}^{\R}_n$  the Gromov boundary has exactly one point and the hyperbolic unfolding roughly looks like a cylinder when we approach infinity. Again the totally geodesic case is trivial: the Gromov boundary is empty since the unfolding is the compact one-point space even when $H$ was a Euclidean hyperplane.
\end{remark}
%
%%%%%%%%%%%%%%%%%%%%%%%%%%%%%%%%%
\subsubsection{Naturality}\label{nat}
%%%%%%%%%%%%%%%%%%%%%%%%%%%%%%%%%
To conclude the introduction we discuss in detail the naturality property (S4) from Definition~\ref{def1} for $H_i\in{\cal{H}}_n$  informally saying  $\bp$ continuously depends on deformations of the underlying space. The extension to almost minimizers is explained in Appendix A.II, III. \\

First, we discuss our notion of convergence for the underlying spaces. Consider a sequence of area minimizers $H_i\in{\cal{H}}_n$ inside a sequence of ambient complete Riemannian manifolds $M_i=M_i^{n+1}$. (We drop any reference to their metrics to ease notation.) We fix base points $p_i\in H_i$. To say that the pointed sequence $(H_i\subset M_i,p_i)$ \emph{converges to the pointed hypersurface} $H:=H^n\subset M:=M^{n+1}$, $p\in H$, means the following.\\

\noindent\textbf{Ambient Level} \,
The $M_i$ compactly $C^k$-converge to a limit manifold $M$ so that $p_i \ra p \in M$ for $i \ra \infty$. This means that for any given $R>0$, there are diffeomorphisms $\Psi_i : B_R(p_i) \ra B_R(p)$ for $i$ sufficiently large, so that $\|\Psi_{i*}g_{M_i} - g_{M}\|_{C^k} \ra 0$ on $B_R(p)$.  In order to have a generous amount of regularity  we generally assume that $k \ge 5$.\\

\noindent\textbf{Minimizer Level} \,
The $H_i$ subconverge to $H$, that is, there exists a convergent subsequence to the limit area minimizer $H \subset M$. This means that for any $R>0$ the sequence $\Psi_i (B_R(p_i) \cap H_i) \ra B_R(p) \cap H$ is subconvergent in $M$ with respect to the  \emph{flat norm topology} (see Equation~\eqref{fll} in Appendix A) one can interpret as a measure for the volume between
$\Psi_i (B_R(p_i) \cap H_i)$ and $B_R(p) \cap H$ within $B_R(p) \subset M$.

\begin{example}\label{ex1} \,
Consider an initial hypersurface $H_0$ in $M_0$ and the rescaled sequence $H_i:=\tau_i \cdot H_0\subset M_i:=\tau_i\cdot M_0$ of \emph{blow-ups} for some $\tau_i \ra \infty$. Fix a singular point $p_0\in\Sigma\subset H_0$ and set $p_i:=p_0$. Then $M_i$ converges compactly to $\R^{n+1}$ and we find a local flat norm subconvergence of $H_i$ to a limit space $H \subset \R^{n+1}$. This limit is actually an area minimizing cone, a so-called \emph{tangent cone}, cf. Appendix \ref{flat-norm-approx}.
\end{example}

\begin{remark} \,
When we do \emph{not} fix the base point, the subconvergence under blow-ups still leads to a limit hypersurface $H \subset \R^{n+1}$. Again, this is a complete area minimizer and oriented boundary in $\R^{n+1}$, but not necessarily a cone.
\end{remark}

Next assume that $B_R(p) \cap H$ is \emph{smooth}. By standard regularity, flat norm convergence of $\Psi_i(B_R(p_i) \cap H_i)$ to $B_r(p)\cap H$ implies that $B_R(p_i) \cap H_i$ is also smooth for sufficiently large $i$ (possibly upon shrinking the radius). Further, flat norm convergence induces $C^k$-convergence in the following sense. Let $\nu\to B_R(p) \cap H$ denote the normal bundle of $B_R(p)\cap H$. Then, for $i$ large enough, the $\Psi_i (B_R(p_i) \cap H_i)$ can be identified with local $C^k$-sections of $\nu$ (cf.\ also Section II in Appendix A and references quoted there). Hence we get canonical diffeomorphisms $\Gamma_i :B_R(p) \cap H \ra \Psi_i (B_R(p_i) \cap H_i)$, and the flat norm convergence implies $C^k$-convergence of the sections $\Gamma_i$ to the \emph{zero section} $B_R(p) \cap H$.

\begin{definition}[$\D$\textbf{-map}]\label{idmap} \,
For sufficiently large $i$, we call the uniquely determined section of $\nu$,
\[
\D := \Gamma_i : B_R(p) \cap H \ra  \Psi_i(B_R(p_i) \cap H_i)
\]
the \textbf{asymptotic identification} or $\D$\textbf{-map} for short.
\end{definition}

\begin{remark}\label{te} \,
1.\ Using finite ball covers we extend the notion of $\D$ maps to domains in $H \setminus \Sigma_{H}$ with compact closures. Writing $id_{H}$ for the zero section of the normal bundle $\nu$, local $C^k$-convergence of the $H_i$ can be rephrased as local $C^k$-convergence of maps, i.e., $|\D - id_{H}|_{C^k} \ra 0$.\\
2.\ In general, $\D(\p B_R(p)) \neq \p B_R(p_i)$, but $\D(\p B_R(p))$ gradually approaches $\p B_R(p_i)$, for $i \ra \infty$.   However, it is only the portion away from these boundaries we are interested in. Thus, whenever needed, we may easily adjust the definitions near the boundary and henceforth ignore these negligible adjustments.
\end{remark}

\begin{definition}[Natural assignements]\label{cr} \,
Consider an assignment $F:H \mapsto F_H$, $H \in {\cal{H}}$, of functions $F_H:H \setminus\Sigma_H\to\R$. Then $F$ is called \textbf{natural}, if $F_H$ commutes with the convergence of underlying spaces. That is, for any pointed sequence $H_i\in {\cal{H}}$, $p_i\in H_i\setminus\Sigma_{H_i}$, locally converging in flat norm to the pointed space $H$, $p\in H \setminus \Sigma_H$ as above, there is a neighborhood $U(p) \subset H \setminus \Sigma_{H}$ such that
\[
|\D^*F_{H_i} - F_{H}|_{C^k(U(p))}=|F_{H_i} \circ \D - F_{H}|_{C^k(U(p))} \ra 0 \mm{ as } i \ra \infty
\]
for some $k=k(F) >0$.
\end{definition}

Thinking of natural assignments as being ``continuous'' with respect to flat norm convergence we can consider more general assignments $H \mapsto F_H$, like tensors or operators, whenever this makes sense. Formally, this can be accomplished through a representation of assigned entity by a set of local coefficient functions.
\begin{example}\label{ex2} \,For instance, we can pull-back via $\D$ intrinsic curvature notions like sectional, Ricci and scalar curvature of $H$, written $sec_H$, $Ric_H$ and $S_H$ respectively. Similarly, the Riemann and Weyl tensors $Riem_H$ and $W_H$ are natural, as well as the extrinsic second fundamental form $A_H$ and its norm $|A_H|$. The Laplace operator $\Delta_H$ and thus the conformal Laplacian $L_H = -\Delta_H +\frac{n-2}{4 (n-1)} \cdot S_H$ and the Jacobi field operator $J_H=-\Delta_H - |A|^2-Ric_M(n,n)$ (where $n$ is the normal vector field of $H$), are also natural assignements. In turn, $dist_{g_H}(\cdot,\Sigma_H)$ is \emph{not} a natural assignment, since the singular set $\Sigma_H$ of the limit $H$ of a converging sequence $H_i$ in ${\cal{H}}$ may have a  different structure than the singular sets $\Sigma_i$ of $H_i$ (consider e.g. a sequence of smooth $H_i$ converging to a singular $H$).
\end{example}
%
%
%
%
%
%%%%%%%%%%%%%%%%%%%%%%%%%%%%%%%%%
%%%%%%%%%%%%%%%%%%%%%%%%%%%%%%%%%
%%%%%%%%%%%%%%%%%%%%%%%%%%%%%%%%%
\setcounter{section}{2}
\renewcommand{\thesubsection}{\thesection}
\subsection{\si-transforms and \si-uniformity}\label{siu}
%%%%%%%%%%%%%%%%%%%%%%%%%%%%%%%%%
%%%%%%%%%%%%%%%%%%%%%%%%%%%%%%%%%
%%%%%%%%%%%%%%%%%%%%%%%%%%%%%%%%%
In this section we construct a concrete family of \si-transforms and establish the \emph{\si-uniformity} for the open manifold $H\setminus\Sigma$.
%
%
%
%%%%%%%%%%%%%%%%%%%%%%%%%%%%%%%%%
\subsubsection{Connectedness of $H\setminus\Sigma$}\label{rn}
%%%%%%%%%%%%%%%%%%%%%%%%%%%%%%%%%
Consider an area minimizing hypersurface $H\in{\cal H}_n$. According to our convention it is connected. Here we want to prove that $H \setminus \Sigma$ is rectifiably connected. Although the codimension of $\Sigma$ in $H$ is greater or than $7$ (see Proposition~\ref{arm} in Appendix A), this is not evident since $H$ degenerates towards $\Sigma$.

\begin{proposition}[Connectedness  of $\mathbf{H \setminus \Sigma}$]\label{coh} \,
For any $H\in{\cal H}_n$ with singular set $\Sigma_H$, the regular complement $H \setminus \Sigma_H$ is rectifiably connected.
\end{proposition}

\noindent{\bf Proof} \,
A connected Riemannian manifold is path connected. Since any continuous curve can be approximated by a rectifiable one, it is sufficient to show that $(H\setminus\Sigma,g_H)$ is connected. For $\Sigma\v$ this is trivial. So let us assume that $\Sigma\n$ and that $H \setminus \Sigma$ contains at least two open, non-empty and disjoint components $C$, $D\subset H\setminus\Sigma$ with $C\cup D=H\setminus\Sigma$. The idea is to think of $C$ and $D$ as minimal currents with boundary $\p C$, $\p D\subset \Sigma$, and to derive a contradiction to the isoperimetric inequality. Towards that end we want to use the local decomposition of a rectifiable current into a locally disjoint collection of oriented minimal boundaries. Concretely, for any $p \in \Sigma$, there is an $r_p >0$ so that $B_{r_p}(p) \cap H \subset M$ is an oriented boundary in $B_{r_p}(p)\subset M$ (cf.\ Section V Appendix A, in particular Proposition~\ref{dicc}). Since this is not a global decomposition we prove the following stronger claim for the case where $H$ is an oriented boundary. \\

\noindent\textbf{Local connectedness} \,
\emph{For any $p \in \Sigma$ and $r \in (0,r_p)$, we choose the connected component $H_r(p)\subset B_r(p) \cap H$  containing $p$. Then $H_r(p)\setminus\Sigma$ is still connected for $r>0$ small  enough.}\\

So assume to the contrary that for arbitrarily small $r \in (0,r_p)$, we have a decomposition $H_r(p)\setminus \Sigma=C_r(p)\cup D_r(p)$ into two open, non-empty and disjoint subsets. (The case of more than two such components can be treated similarly.) We rescale $B_r(p)\subset M^{n+1}$ to unit size so that $B_1:=B_1(p) \subset r^{-1} \cdot M$ approximates a Euclidean $(n+1)$-ball as closely as we wish, let us say as in Ch.\ref{nat}, in $C^5$-topology.  We denote by $H_1$ and $C_1$, $D_1$ the rescaled oriented boundary $H_r(p)$ and components $C_r(p)$ and $D_r(p)$. Inside $B_1$ we can choose a tubular neighborhood $\U$ of $\Sigma$
such that $\p\U\cap D_1$ is smooth with $Vol_n(\p\U\cap D_1)\ra0$ if we shrink $\U$ towards $\Sigma$. (One may define such $\U$ as distance tubes of some mollified distance function. Then one uses the coarea formula \ \cite[2.1.5,Theorem 3]{GMS} and that the Hausdorff dimension of $\Sigma$ is $\le n-7\le n-2$.)\\

Next let $C^*_1:= C_1\cup(\U\cap H_1)$ and $D^*_1:= (H_1\setminus \Sigma)\setminus C^*_1$. Then both $C^*_1$ and $D^*_1$ are integral currents, and we have the decomposition $H = C^*_1 \cup D^*_1$ with $Vol_n(H) = Vol_n(C^*_1) + Vol_n(D^*_1)$. Now a variant of the \emph{isoperimetric inequality} for oriented minimal boundaries due to Bombieri and Giusti~\cite[Thm.\ 2, p.\ 31]{BG} gives
\[
Vol_{n-1}(\p D^*_1 \cap B_1(p)) \ge  k_n \cdot \min \{Vol_n(C^*_1 \cap  B_{\beta_n}(p)), Vol_n(D^*_1 \cap  B_{\beta_n}(p))\}^{(n-1)/n},
\]
for some constants $k_n > 0$ and $\beta_n \in (0,1)$ depending only on the dimension $n$. (In [BG] this inequality is formulated for minimal boundaries in thee Euclidean space, but the argument carries over to minimal boundaries in $B_1(p) \subset r^{-1} \cdot M$ for $r$ large enough.) In particular, the right hand side is positively lower bounded which contradicts $Vol_{n-1}(\p\U\cap D_1) \ra 0$ if we shrink $\U$ towards $\Sigma$.
\qed

%
%
%
%%%%%%%%%%%%%%%%%%%%%%%%%%%%%%%%%
\subsubsection{\si-transforms}\label{st}
%%%%%%%%%%%%%%%%%%%%%%%%%%%%%%%%%
The easiest way to define \si-transforms is to use distance tubes of the $|A|$-level sets on $H$. An alternative approach, which we will not discuss here, is to choose area minimizing hypersurfaces \emph{within} $H \setminus \Sigma$ which are spanned over the \emph{obstacle} $|A|^{-1}[c,\infty)$.\\

\noindent\textbf{Metric \si-transforms} \,
Choose $\alpha >0$, $c>0$, and let $H \in {\cal{H}}$ be a non-totally geodesic area minimizer. We first define the $|A|$-\textbf{skins} $\M_c=\M_c(\alpha)$ of the desired \si-transform $\bp_{\alpha}$ by
\[
\M_c(\alpha):= \mm{\emph{ the boundary of the} } \alpha/c \mm{\emph{-distance tube }} \U^\alpha_c \mm{ \emph{of }} |A|^{-1}[c,\infty)
\]
where the distances are measured with respect to $d_{g_H}$.\\

For $c<d$ we have $\overline{\U}^\alpha_d\subset \U^\alpha_c$ and therefore $\M_c \cap \M_d \v$, since $|A|^{-1}[d,\infty)\subset |A|^{-1}[c,\infty)$ and $\alpha/d < \alpha/c$. We can thus uniquely define
\[
\bp_{\alpha,H}(x):= c\mm{ for } x\in \M_c.
\]
Usually $\bigcup_{c > 0} \M_c \subsetneq H \setminus \Sigma$ but  the definition can be canonically extended as follows.

\begin{definition}[Metric \si-transforms]\label{lsk} \,
For $\alpha > 0$ and $H\in\cal{H}$ we define the \textbf{metric \si-transform} $\bp_{\alpha,H}$ as follows. When $H$ is totally geodesic, we set $\bp_{\alpha,H}\equiv 0$. Otherwise, we let
\[
\bp_{\alpha,H}(x):= \sup \{c \,| \, x \in \overline{\U}^\alpha_c\}
\]
for any $x \in H \setminus \Sigma$. In order to ease notation we usually write $\bp_\alpha$ if the associated area minimizer $H$ is clear from the context.
\end{definition}

\begin{lemma}[Divergence of  $\bp_{\alpha}$]\label{diva} \, For any  sequence $p_i \in H \setminus \Sigma_H$ and $p \in \Sigma_H$ with
\[d_{g_H}(p_i,p)= dist_{g_H}(p_i,\Sigma_H) \ra 0\mm{ for } i \ra \infty,\mm{ we have } \bp_{\alpha}(p_i) \ra \infty.\]
\end{lemma}

\noindent{\bf Proof} \,  We assume we had a converging sequence of points $p_i \in H \setminus \Sigma$ and some limit $p \in \Sigma$ with $d_{g_H}(p_i,p)= dist_{g_H}(p_i,\Sigma) \ra 0$ and  $\bp_{\alpha}(p_i) <c$, for some common $c>0$, when  $i \ra \infty$. That is, we have $dist_{g_H}(p_i, |A|^{-1}[c,\infty)) > \alpha/c$, for all  $i$. For $2 \cdot d_{g_H}(p_i,p) <  \alpha/c$ we infer that $|A| <  c$  on $B_{2 \cdot d_{g_H}(p_i,p)}(p_i) \setminus \Sigma_H$ and, hence, on $B_{d_{g_H}(p_i,p)}(p)\setminus \Sigma_H$. This contradicts the assumption $p \in \Sigma_H$, since  \ref{t3} shows that $|A|$ is unbounded near singular points.\qed

\begin{definition}[\si-distance] \,
For $H$ non-totally geodesic we define the \textbf{\si-distance} by
\[
\delta_{\bp_{\alpha}} := 1/\bp_{\alpha}:H\setminus\Sigma\to\R.
\]
For $H$ totally geodesic we set accordingly $\delta_{\bp_{\alpha}}\equiv\infty$ (cf.\ Definition~\ref{def1}).
\end{definition}

\begin{proposition}[Relations between $\mathbf{\bp_\alpha, |A|}$ and distance functions]\label{locskin} \,
For any non-totally geodesic $H$ we have the following estimates on $H \setminus \Sigma$.
\begin{description}
  \item[A.] \emph{\textbf{Growth estimates and Lipschitz properties of the \si-distance}}\, The \si-distance $\delta_{\bp_{\alpha}}$ is $1/\alpha$-Lipschitz on $H\setminus\Sigma$:
  \begin{equation}\label{li}
  |\delta_{\bp_{\alpha}}(p)- \delta_{\bp_{\alpha}}(q)|   \le d_{g_H}(p,q)/\alpha, \mm{ in particular } \delta_{\bp_{\alpha}}(x) \le dist_{g_H}(x,\Sigma)/\alpha.
  \end{equation}

  For totally geodesic $H$, we set  $|\delta_{\bp_{\alpha}}(p)- \delta_{\bp_{\alpha}}(q)|:=0$  to make (\ref{li}) consistent on ${\cal{H}}$.
  \item[B.] \emph{\textbf{Interpolation properties of the metric \si-transforms}}
  \begin{enumerate}
    \item $\bp_{\alpha} \ge \bp_{\beta}, \mm{ for }\, \alpha \ge \beta >0$.
    \item $\bp_{\alpha} \ra |A|$ in $L_{loc}^\infty$ as $\alpha \ra 0$.
    \item $\bp_{\alpha}/\alpha  \ra 1/ dist_{g_H}(\cdot,\Sigma)\mm{ in } L_{loc}^\infty, \mm{ as }\,\alpha \ra \infty$.
  \end{enumerate}
\end{description}
\end{proposition}

\noindent{\bf Proof} \textbf{A.} \,
We may assume that $p \in \M_c$ and $q \in \M_d$ for some $d > c >0$. Then
\[
|\delta_{\bp_{\alpha}}(p)- \delta_{\bp_{\alpha}}(q)| = |1/\bp_{\alpha}(p)- 1/\bp_{\alpha}(q)| =  \alpha^{-1} \cdot \left|\frac{\alpha}{c} - \frac{\alpha}{d}\right| \le d_{g_H}(p,q)/\alpha.
\]
The latter inequality follows from $\overline{\U}^\alpha_d\subset \U^\alpha_c$. From this we also infer  $\delta_{\bp_{\alpha}}(p) \le dist_{g_H}(p,\Sigma)/\alpha$ by moving $q$ towards $\Sigma$ from \ref{diva}.\\

\noindent\textbf{B.} \,
The inequality (i) follows from $\U^\beta_c \subset \U^\alpha_c$ for $\alpha \ge \beta >0$. The boundary of the $\alpha/c$-distance tube of $|A|^{-1}[c,\infty)$ converges locally uniformly to $|A|^{-1}(c)$ for $\alpha \ra 0$, whence $\bp_{\alpha} \ra |A|$ in $L_{loc}^\infty$ as $\alpha\ra0$. For (iii), we note that $dist_{g_H}(x, |A|^{-1}[d,\infty))\le \alpha/d$ if we set $d:= \bp_\alpha(x)$. Since $|A|^{-1}[d,\infty)$ shrinks to $\Sigma$ as $d \ra \infty$ the claimed convergence follows.
\qed

\begin{proposition}\label{locskintr} \,
$\bp_\alpha$ is an \si-transform for any $\alpha > 0$.
\end{proposition}

\noindent{\bf Proof} \,
We need to verify the axioms (S1) - (S4) from Definition~\ref{def1}.\\

\textbf{(S1)} and \textbf{(S2)}: \,
From the definition $\bp_\alpha \equiv 0$ if $H \subset M$ is totally geodesic, and $\bp_\alpha>0$ with $\bp_\alpha \ge |A|$ if not. Under scalings of $H$ by $\lambda$, the distances on $H$ and the function $|A|$ scale by $\lambda$ and $1/\lambda$, respectively. Hence $\bp_{\alpha, \lambda \cdot H} \equiv \lambda^{-1} \cdot  \bp_{\alpha, H}$. Finally, the divergence of $\bp_\alpha$ towards $\Sigma$ was checked in \ref{diva}\\

\textbf{(S3):} \,
The Lipschitz regularity of $\delta_{\bp_\alpha}:=1/\bp_{\alpha}$ with Lipschitz constant $1/\alpha$ is just Proposition~\ref{locskin}~(A).\\

\textbf{(S4):} \,
The naturality of $\bp_{\alpha}$ follows from that of $|A|$ and standard regularity theory. If for $H_i$, $H\in {\cal{H}}$ we express convergence of pointed spaces $(H_i,p_i)\ra (H,p)$ on some small ball in $H\setminus\Sigma$ in terms of $\D$-maps, we obtain compact smooth convergence on $H \setminus \Sigma_{H}$ in virtue of Remark~\ref{te} (i). The naturality of $|A|$ yields compact convergence of $|A_{H_i}|$ and thus compact $L^{\infty}$-convergence of $\bp_{\alpha,H_i}$. Indeed, there are two cases to consider. Either $\bp_{\alpha}(p_i) \ra 0$. But then the $H_i$ converge to a totally geodesic limit with $\bp_{\alpha,H}\equiv0$. Or there exists a converging sequence $p_i \in H_i \setminus \Sigma_{H_i}$ so that $\bp_{\alpha}(p_i)$ remains positively bounded from below by some $d >0$. Then $p_i \in\U^\alpha_d$ and thus $p \in \U^\alpha_d \subset H$. Hence convergence of $|A|$ implies convergence of the values $\bp_{\alpha}(p_i)$. Since we have a uniform Lipschitz estimate for $\delta_{\bp_\alpha,H}$, $H \in {\cal{H}}$, Rellich compactness yields $C^\gamma$-H\"older subconvergence, for any $\gamma \in (0,1)$. This is actually convergence for we have a well-defined $L^\infty$-limit.
\qed

The latter result give us a working model of an \si-transform very much as singular homology shows that there at least one theory that satisfies the Eilenberg-Steenrod axioms. From this point on, we no longer refer to any particular model $\bp$ of an \si-transform but  derive all further results exclusively from the axioms (S1)-(S4).
%
%
%%%%%%%%%%%%%%%%%%%%%%%%%%%%%%%%%
\subsubsection{Uniformity and \si-Uniformity}\label{inta}
%%%%%%%%%%%%%%%%%%%%%%%%%%%%%%%%%

To control the geometry of $H \setminus \Sigma$ near $\Sigma$ the language of uniform spaces turns out to be very natural. Consider a non-complete, locally compact and complete, rectifiably connected metric space $(X,d_X)$. For such a space we let $\overline{X}$ denote its metric completion and set $\p X:= \overline{X} \setminus X$. Recall from Definition~\ref{ud} that $X$ is called a \emph{c-uniform space}, or \emph{uniform space} for short, if there exists a constant $c \ge 1$ such that any two points can be joined by a \emph{c-uniform curve} in $X$. This is a rectifiable curve $\gamma: [a,b] \ra X$, for some $a<b$, running from $p$ to $q$ such that $\gamma$ is \emph{quasi-geodesic} satisfying the \emph{twisted double cones condition}, i.e.,
\[
l_{d_X}(\gamma)  \le c \cdot  d_X(p,q) \mm{ and } l_{min}(\gamma_{p,q}(z)) \le c \cdot dist_{d_X}(z,\p X)\mm{ for any }z \in \gamma_{p,q},
\]
where $l_{min}(\gamma_{p,q}(z))$ is the minimum of the lengths of the two subcurves of $\gamma_{p,q}$ from $p$ to $z$ and from $q$ to $z$.\\

\noindent\textbf{\si-uniform spaces.} \, We consider the space $X = H \setminus \Sigma$ with $\p X= \Sigma$ and ask if $X$ is a uniform space. This is actually the case, but we prove this claim we need to go still one step further and first establish an \emph{\si-uniformity} of $X$.\\

To better understand this strategy, we observe that complexity and curvature of $H \setminus \Sigma$ are not properly coupled to the metric distance to $\Sigma$.
Therefore it seems rather delicate to approach the proof of the desired metric twisted cone condition
\[l_{min}(\gamma_{p,q}(z)) \le a \cdot dist_{g_H}(z,\p X)\] directly. Instead, we employ the naturality of \si-structures to use compactness results for ${\cal{H}}^n$. They allows us to derive a sharpened \si-version for a given  \si-transform $\bp$:
\[l_{min}(\gamma_{p,q}(z)) \le b \cdot \delta_{\bp}(z),\] and, a posteriori, we infer the result also for $dist_{g_H}(z,\p X)$, from the general relation $\delta_{\bp}(z) \le c \cdot dist_{g_H}(z,\p X)$. We point out that using an \si-version of the twisted cone condition is more than a technicality. The hyperbolic unfoldings (and the analytic applications mentioned in the introduction) rely on this stronger \si-version of uniformity.\\

To formulate our main result we fix some \si-transform $\bp$. Also for our notational convenience we assume that the Lipschitz constant for $\delta_{\bp}$ equals $1$.

\begin{proposition}[\si-uniformity of $\mathbf{H \setminus \Sigma}$]\label{intsuni} \,
For any connected hypersurface $H \in {\cal{H}}^n$ with (possibly empty) singular set $\Sigma=\Sigma_H$, we have
\begin{enumerate}
  \item $H \setminus \Sigma$ and $H$ are rectifiably connected. In particular, any compact $H \in {\cal{H}}$ has a finite intrinsic diameter: $diam_{g_H}H < \infty$.
  \item For some $c >0$, $H \setminus \Sigma$ is a c-\si-uniform space. That is, any pair $p,q \in H \setminus \Sigma$ can be joined by a c-\si-uniform curve in $H \setminus \Sigma$, i.e., a rectifiable curve $\gamma_{p,q}: [a,b] \ra H \setminus \Sigma$, for some $a <b$, with $\gamma_{p,q}(a)=p$, $\gamma_{p,q}(b)=q$, so that the following conditions hold:
  \begin{itemize}
    \item \emph{\textbf{Quasi-geodesic}:} \,
    $l_{g_H}(\gamma)  \le c \cdot  d_{g_H}(p,q)$.
    \item \emph{\textbf{Twisted double \si-cones}:} \,
    $l_{min}(\gamma_{p,q}(z)) \le c \cdot \delta_{\bp}(z)$ for any $z \in \gamma_{p,q}$.
  \end{itemize}
  \item More generally, any pair $p$, $q \in H$ can be joined by a c-\si-uniform curve supported in $H \setminus \Sigma$, except for its endpoints if $p$ or $q\in\Sigma$.
 \item For $H \in {\cal{H}}^{\R}_n$ we get a uniform constant $c_n$ depending only on the dimension $n$ so that $H \setminus \Sigma$  is a $c_n$-\si-uniform space.
\end{enumerate}
\end{proposition}

From (S3) in Definition~\ref{def1} we immediately draw the

\begin{corollary}[Uniformity of $\mathbf{H \setminus \Sigma}$]\label{intuni} \,
For any singular area minimizer $H \in {\cal{H}}^n$  the metric space $H \setminus \Sigma$ is uniform.
\end{corollary}

\begin{remark}
The distance function $dist_{g_H}(\cdot,\Sigma)$ does not behave naturally under convergence of the underlying spaces so that our subsequent strategy in the \si-uniform setting does not apply directly to the uniform setting.
\end{remark}

These \emph{intrinsic} uniformity properties of $H \setminus \Sigma$ clearly rely on the \emph{extrinsic} property of $H$ to be an area minimizer in its ambient space. We therefore also get the extrinsic estimates when $H$ is an oriented boundary.

\begin{corollary}[Intrinsic versus extrinsic metric]\label{fini} \,
For any $H\in{\cal{H}}^{\R}_n$ and $p$, $q\in H\subset\R^{n+1}$ there is a constant $c^{\R}_n \in (0,1)$ depending only on the dimension $n$ and such that
\[
c^{\R}_n \cdot  d_{g_H}(p,q) \le  d_{g_{\R^{n+1}}}(p,q) \le d_{g_H}(p,q).
\]
\end{corollary}

Finally, hypersurfaces pass their c-uniformity constant to their blow-up limits.

\begin{corollary}[Inheritance under blow-ups]\label{insu} \,
For $H\in {\cal{H}}^n$ let $H\setminus \Sigma$ be a $c$-\si-uniform space for some $c >0$. If $F$ is any blow-up limit of $H$, then $F \setminus \Sigma_F$ is also $c$-\si-uniform.
\end{corollary}
%
%
%
%%%%%%%%%%%%%%%%%%%%%%%%%%%%%%%%%
\subsubsection{From isoperimetry to quasi-geodesic pipelines}\label{intapr}
%%%%%%%%%%%%%%%%%%%%%%%%%%%%%%%%%
Let us now start with the proof of Proposition~\ref{intsuni} as well as its corollaries.
To build \si-uniform curves we shall proceed in several steps. We gradually upgrade the rectifiable connectedness of $H \setminus \Sigma$, cf. Proposition~\ref{coh}, until we reach the asserted \si-uniformity. We use the \emph{scaling invariance} of the area minimizing condition, and the \emph{naturality} of $\bp$. To avoid trivialities we assume to work with non-totally geodesic hypersurfaces.\\

\noindent\textbf{Step 1 (Short quasi-geodesic curves)}\label{co1}\\
Here we derive the existence of \emph{short quasi-geodesic curves} with some \emph{controlled \si-distance} to the boundary. We start with hypersurfaces in $\R^{n+1}$. This corresponds to the limit case of strong rescalings of $H \subset M$. For $\rho>0$ set
\[
\I(\rho):= \{x \in H \setminus \Sigma \,|\, \delta_{\bp}(x) < \rho \}\quad\mm{and}\quad\E(\rho):= H \setminus \I(\rho).
\]
We think of a point $p\in\E(\rho)$ as being at ``\emph{curved distance}'' $\delta_{\bp}(p)\geq\rho$ from $\Sigma$. Unlike the metric distance it is stable under perturbations of $H$ within ${\cal{H}}_n$ via $\D$-maps.

\begin{lemma}\label{copr} \,
Let $H\in{\cal{H}}^\R_n$. Then for any $t >0$ there are some $\tau\equiv\tau(t,n) \in (0,t)$ and $\Pi\equiv \Pi(t,n)>0$ such that any two points $p$, $q\in\E(t)\subset H$ with $d_{g_{\R^{n+1}}}(p,q) =1$ can be connected by a \emph{\textbf{short quasi-geodesic curve}} $\gamma_{p,q}\subset\E({\tau})$ of length $l_{g_H}(\gamma_{p,q}) \le \Pi$.
\end{lemma}

\noindent\textbf{Proof} \,
Assume that there are \emph{no} estimates as in Lemma~\ref{copr} which are valid for all hypersurfaces. Then we can pick a sequence of such hypersurfaces $H_i\subset \R^{n+1}$, bounding connected open sets $U_{H_i}\subset \R^{n+1}$ \cite[Theorem 1 p.\ 26 and Corollary p.\ 30]{BG} as well as points $p_i$, $q_i \in  \E(t) \subset H_i$ with $d_{g_{\R^{n+1}}}(p_i,q_i) =1$, such that the intrinsic distance $d_{\E({1/i})}(p_i,q_i)$ in $\E({1/i})$ diverges: Either $d_{\E(1/i)}(p_i,q_i) \ge i$ or $=\infty$ if $\E({1/i})$ is not connected and $p_i$, $q_i$ lie in different components. We may assume that $p_i=0= (0,\ldots,0)$, $q_i = e_1 =(1,0,\ldots,0)$, and that $\{H_i\}$ converges compactly to a limit area minimizer $H_\infty$~\cite[Theorem 1.19 and Lemma 9.1]{Gi}. The Lipschitz estimate $|\delta_{\bp}(p)- \delta_{\bp}(q)|\le L \cdot d_{g_{H_i}}(p,q)$ shows that $B_r(p_i)$, $B_r(q_i) \subset \E(t/2)$ for any $r \in (0, t/(2 \cdot L))$. The non-extinction statement from Proposition~\ref{nex} implies that these balls are not annihilated in the limit so that the limit points of $\{p_i\}$ and $\{q_i\}$, namely $0$ and $e_1$, belong to $\E(t) \subset H_\infty$. Now $H_\infty \setminus \Sigma_{H_\infty}$ is rectifiably connected and $\bp$ is a proper function on $\overline{B}_R(0) \cap H_\infty \setminus\Sigma_{H_\infty}$ for any given $R >0$, since $\delta_{\bp}(x) \le L  \cdot dist_{g_{H_\infty}}(x,\Sigma_{H_\infty})$. Hence, there is a smooth curve $\gamma$ in $H_\infty$ which connects $0$ and $e_1$ within $\E(\tau) \subset H_\infty$ for some suitably small $\tau \in (0, t)$. Again the properness of $\bp$ on $\overline{B}_R(0) \cap H_\infty \setminus\Sigma_{H_\infty}$ shows that there is a smooth tube $\U\subset \E(\tau/2)$ around $\gamma$. We thus get smooth convergence of suitable tubes $\U_i \subset H_i$ to $\U$ and rectifiable curves $\gamma_i \subset\U_i$, connecting $p_i$ and $q_i$, to the curve $\gamma\subset H_\infty$. For large $i$, the naturality of $\bp$, axiom (S4),  shows that $\U_i  \subset \E(\tau/4) \subset H_i$ and $l_{g_{H_i}}(\gamma_i) \le l_{g_{H_\infty}}(\gamma)+1$, contradicting the assumption.
\qed

\noindent\textbf{Step 2 (Pipelines of short quasi-geodesics curves)}\\
Next we assemble the short quasi-geodesics to form quasi-geodesic ``\emph{pipelines}'' in sufficiently small but uniformly sized balls in $H$. We explicitly allow these balls to be centered in $\Sigma_H$. We first establish a basic volume control.

\begin{lemma}\label{coprq} \,
Let $H\in{\cal{H}}^\R_n$. Then for any $\varpi\in (0,1)$ there is some $t_{n,\varpi} \in (0,1)$ such that for any $t \in (0, t_{n,\varpi})$, $k \in \Z$ and $p$, $q\in H$ with $d_{\R^{n+1}}(p,q)=3/2$, we have
\begin{align}
0<\varpi&\le\frac{Vol_n(\E(2^{-k} \cdot t) \cap B_{2^{-k+1}} \setminus B_{2^{-k}}(p)\cap H)}{Vol_n(B_{2^{-k+1}}\setminus B_{2^{-k}}(p)\cap H)}\label{w1}\\[5pt]
0<\varpi&\le\frac{Vol_n(\E(t/2) \cap B_{1}(p) \cap B_{1}(q) \cap H)}{Vol_n(B_{1}(p) \cap B_{1}(q) \cap H)}\label{w2}
\end{align}
where $B_R\setminus B_r(z):=B_R(z)\setminus B_r(z)$ is the difference of the balls of radius $R$ and $r$ in $\R^{n+1}$.
\end{lemma}

In particular, for any given point $z\in B_{2^{-k+1}} \setminus B_{2^{-k}}(p)$, any $\ve >0$ and $\varpi(\ve)$ sufficiently close to $1$, we have $\E(2^{-k} \cdot t) \cap B_{2^{-k} \cdot \ve}(z)\not=\emptyset$ for $t \in (0, t_{n,\varpi})$.\\

\noindent\textbf{Proof} \,
The isoperimetric inequality~\cite[5.13, 5.14 and Inequality (5.16)]{Gi} and a simple comparison with the (larger) volume of $\p B_1(0)$ give positive constants $c_n^\pm$ depending only on the dimension $n$, such that for any $p\in H$ we have $c_n^- \le Vol_n(B_{1/2}(p) \cap H) \le c_n^+$. On the other hand, we have for any given $H$ and $p\in H$ that
\[
\frac{Vol_n(\E(t/2) \cap B_{1}  \setminus B_{1/2}(p) \cap H)}{Vol_n(B_{1} \setminus B_{1/2}(p) \cap H)} \ra 1\mm{ for } t \ra 0 \mm{ since }\bigcap_{t >0}\I(t/2) = \Sigma \mm{ and } Vol_n(\Sigma) = 0.
\]
Compactness arguments for area minimizers analogously to those in Lemma~\ref{copr} then yield some $t_{n,\varpi} \in (0,1)$ such that for any $t \in (0, t_{n,\varpi})$ and for any point $x \in H$
\[
Vol_n(\E(t/2) \cap B_{1}  \setminus B_{1/2}(x) \cap H) \ge \varpi \cdot Vol_n(B_{1} \setminus B_{1/2}(x) \cap H) \ge \varpi \cdot c_n^- .
\]
This implies~\eqref{w1} for $k=0$. The case $k \neq 0$ follows from scaling by $2^k$ and the scaling behaviour of $\bp$ and the volumes. Inequality~\eqref{w2} can be derived in the same way.
\qed

From now on we consider a general hypersurface $H \in {\cal{H}}^n$. In particular, $H\subset M$ could be compact. We assume, as we already did on several occasions, that $M$ has been scaled by some large constant to the effect that for some $\varpi$ sufficiently close to $1$, every ball $B_r\subset M$ of radius $r \le 10 \cdot \Pi(t_{n,\varpi},n) + 5$ is very close to the ball $B_r(0) \subset \R^{n+1}$ in some sufficiently regular topology e.g., in  $C^5$-topology. In particular, we may apply Lemmas \ref{copr} and \ref{coprq} to $H\cap B_r$ which again we think of as an oriented boundary inside $\R^{n+1}$.\\

Now we explain how to join any given point $p$ with points in $B_1 \setminus B_{1/2}(p) \cap H$ by a quasi-geodesic in $H$ with controlled \si-distance, i.e., a curve surrounded by a twisted \si-cone pointing to $p$.\\

 From Inequality~\eqref{w1} we may assume we have some $p_k \in \E(2^{-k} \cdot t_{n,\varpi}) \cap B_{2^{-k+1}} \setminus B_{2^{-k}}(p) \cap H$, for any $k \ge 0$, with $d_{g_M}(p_k,p_{k+1}) = 2^{-k}$. (Indeed this holds up to a multiple arbitrarily close to $1$ and common for all $k$.  ~\ref{coprq} shows that for any $\eta \in (0,1)$  there is a $\varpi$ sufficiently close to $1$ so that the $p_k$ can be chosen in such way that $(1 - \eta) \cdot  2^{-k} \le d_{g_M}(p_k,p_{k+1}) \le (1 + \eta) \cdot  2^{-k}$.)  From Lemma~\ref{copr} we get for some $\tau(t_{n,\varpi},n) < t_{n,\varpi}$ a curve $\gamma_{p_k,p_{k+1}}$ connecting $p_k$ and $p_{k+1}$ with $l_{g_H}(\gamma_{p_k,p_{k+1}}) \le \Pi/2^k$ and $\gamma_{p_k,p_{k+1}} \subset \E(2^{-k} \cdot \tau)\subset H$. Since $\Pi$ may be much larger than $1$ we usually have $\gamma_{p_k,p_{k+1}} \varsubsetneq \E(2^{-k} \cdot \tau) \cap B_{2^{-k+1}} \setminus B_{2^{-k}}(p) \cap H.$ However, writing $\tau^*:=  \tau \cdot 2^{-m}$ with
$m:= \mm{\emph{the smallest integer}}\ge \log_2(10 \cdot \Pi)$, we get
\begin{equation}\label{vk01}
\gamma_{p_k,p_{k+1}} \subset  \E(2^{-(k-m)} \cdot \tau^*) \cap B_{2^{-(k-m)+1}}(p)\cap H,
\end{equation}
since $l_{g_H}(\gamma_{p_k,p_{k+1}}) \le \Pi/2^k$. Now we glue the curves $\gamma_{p_k,p_{k+1}}$, for all $k \ge 0$, by identifying the endpoint of $\gamma_{p_k,p_{k+1}}$ with the starting point of $\gamma_{p_{k-1},p_{k}}$. As a result, we obtain the \textbf{quasi-geodesic pipeline} $\Gamma=\Gamma_{p_0,p}$ from $p_0$ to $p$. Each point $z$ on the subcurve $\gamma_{p_k,p_{k+1}} \subset \Gamma$ remains within a distance $\le \Pi/2^k$ to both endpoints $p_k$ and $p_{k+1}$. In this way we get estimates for any point $z \in \gamma_{p_k,p_{k+1}} \subset\Gamma$ and the subcurve $\Gamma(z) \subset \Gamma$ from $p$ to $z$, namely
\begin{align}
d_{g_H}(p,z) \le l_{g_H}(\Gamma(z))&\le \sum_{a \ge k} \Pi/2^a = \Pi/2^{k-1} \le  \Pi \cdot d_{g_M}(p,z) \le  \Pi \cdot d_{g_H}(p,z)\label{vk001} \\
l_{g_H}(\Gamma(z)) &\le 2 \cdot  \Pi/ \tau \cdot \delta_{\bp}(z).\label{vk1}
\end{align}
For the latter inequality we use \eqref{vk01}. It shows that $\Gamma(z) \subset \E(2^{-k} \cdot \tau)$. Then we get from \eqref{vk001}  $\delta_{\bp}(z) \ge 2^{-k} \cdot \tau \ge \tau/(2 \cdot \Pi) \cdot l_{g_H}(\Gamma(z))$.
\qed

\noindent\textbf{Step 3A (\si-uniformity for Euclidean hypersurfaces)}\\
We use the pipelines of Step 2 to derive the \si-uniformity for $H\in{\cal{H}}^{\R}_n$. Pick $p$ and $q$ in $H$. Since we are in a scaling invariant situation (in particular, the \si-uniformity condition is scaling invariant), we may assume that $d_{g_{\R^{n+1}}}(p,q)=3/2$. From \eqref{w1} we may choose a common starting point $p_0\in B_1(p) \cap B_1(q)\cap H$ to construct pipelines $\Gamma_{p_0,p}$ and $\Gamma_{p_0,q}$. Then \eqref{vk1} shows that the composition of these pipelines defines a $c$-\si-uniform curve from $p$ to $q$ with $c:=4 \cdot \Pi/ \tau + 4 \cdot \Pi$. For a given \si-transform this number depends only on the dimension since we are in the flat Euclidean space and may choose one common $\varpi$ for any point $p \in H$ and also for any $H\in{\cal{H}}^{\R}_n$.\\

Since any two points $p,q \in H$ can be joined by some $c$-\si-uniform curve, supported in $H \setminus \Sigma \cup \{p,q \}$, we find that $H$ is rectifiably connected.
\qed

\noindent\textbf{Step 3B (\si-uniformity for compact hypersurfaces)}\\
For compact hypersurfaces in ${\cal{H}}^c_n$ we combine the pipeline construction with the following consequence of Proposition~\ref{coh}, the connectedness of $H \setminus \Sigma$, to check their \si-uniformity.

\begin{lemma}\label{rnpc} \,
Let $H \in {\cal{H}}^c_n$ and $K \subset H \setminus \Sigma$ be compact. Then there exist $l=l(K) >0$ and $s = s(K) > 0$ such that any two points $p$, $q\in K$ can be linked by a rectifiable curve in $\E(s)$ of length less or equal than $l$.
\end{lemma}

\noindent{\bf Proof} \,
We can cover $K$ by a finite collection of small balls $B_\rho(p_1), ..,B_\rho(p_k)$, $\rho \ll 1$, so that $\overline{B}_{2 \cdot\rho(p_1)},\ldots,\overline{B}_{2 \cdot \rho(p_k)}\subset H \setminus \Sigma$. We link any two centers of these balls, as well as the points $p$ and $q$ to the centers of a ball they belong to, by rectifiable curves in $H \setminus \Sigma$. These curves give $k \, !$ compact sets $\Gamma_1,\ldots,\Gamma_{k \, !}$ in $H \setminus \Sigma$. Since the union of the curves $\Gamma_i$ and the balls $\overline{B}_{2 \cdot \rho(p_i)}$ is compact, we can find some small $s >0$ such that $\overline{B}_{2 \cdot \rho(p_1)}\cup\ldots\overline{B}_{2 \cdot \rho(p_k)}\cup \Gamma_1 \cup\ldots\Gamma_{k \, !} \subset \E(s)$. Then set $l(K) := \max \{l_{g_H}(\Gamma_1),\ldots,l_{g_H}(\Gamma_{k \, !})\}+2$.
\qed

We consider $p$, $q\in H \setminus \Sigma$ and write $\Delta:= d_{g_M}(p,q)$. For a given $\ve>0$ we may scale $M$ in such a way that on $B_5(0) \subset T_pM$ the exponential map $exp_p$ is $\ve$-close to an isometry in $C^5$-norm. We distinguish to cases.\\

\noindent$\mathbf{\Delta \le 1}:$ \,
Choose a local decomposition of $H$ into oriented minimal boundaries such that $p$, $q\in H \setminus \Sigma$ belong to the same boundary. Then we can argue as in the previous step. Upon scaling $M$ by $\frac{3}{2\Delta} \ge 1$ we may assume that $\Delta = 3/2$. Again, \eqref{w1} shows the existence of a common starting point $p_0 \in B_1(p) \cap B_1(q)\cap H$. We can thus construct the pipelines  $\Gamma_{p_0,p}$ and $\Gamma_{p_0,q}$. Using \eqref{vk1} we see that the composition of these two curves is a $c$-\si-uniform curve from $p$ to $q$ with $c:=4 \cdot \Pi/ \tau + 4 \cdot \Pi$.\\

\noindent$\mathbf{\Delta > 1}:$ \,
Let $\U(\Sigma)=\{p\in H\mid dist_{g_H}(p,\Sigma)<1/4\}$ and consider the compact set $K = H \setminus \U(\Sigma)$. If $p$ or $q\in\U(\Sigma)$, we first choose points $p_K$, $q_K \in K$ minimizing the distance to $p$ and $q$. Moreover, we take two local decompositions of $H$ into oriented minimal boundaries such that $p$ and $p_K$ as well as $q$ and $q_K$ belong to the same boundary of one of these decompositions (but not necessarily $p$ and $q$). Now we construct the pipelines $\Gamma_{p_K,p}$ and $\Gamma_{q_K,q}$ with endpoints in $p$ and $q$ as in Step 2. Of course, these curves boild down to constant curves whenever $p$ or $q\notin\U(\Sigma)$. Next, Corollary \ref{rnpc} gives us some $l=l(K) >0$ and $s = s(K) > 0$ so that $p_K$ and $q_K\in K$ can be linked by a rectifiable curve $\gamma_{p_K,q_K}$ in $\E(s)$ of length $\le l$. Since $d_{g_H}(p_K,q_K)>1/2$, this can be rewritten as
{\small \begin{equation}\label{pre}
l_{g_H}(\gamma_{p_K,q_K})\le l = l/d \cdot d \le 2 \cdot l \cdot d_{g_H}(p_K,q_K)  \mm{ and }  l_{min}(\gamma_{p_K,q_K}(z)) \le l \le l/s \cdot s \le l/s \cdot \delta_{\bp}(z).
\end{equation}}
Then we stick the pipelines $\Gamma_{p_K,p}$ and $\Gamma_{q_K,q}$ together with $\gamma_{{x_K,y_K}}$. This defines a curve $\Gamma_{p,q}$ that links $p$ with $q$.\\

We check its \si-uniformity properties: The curve $\Gamma_{p,q}$ is quasi-geodesic as follows from \eqref{vk1}, \eqref{pre} and the inequalities $d_{g_H}(p,p_K)$, $d_{g_H}(q,q_K) \le 1/4 \le d_{g_H}(p,q)/4$:
\begin{align*}
l_{g_H}(\Gamma_{p,q})&\le 2 \cdot \Pi/2 + 2 \cdot l \cdot d_{g_H}(p_K,q_K)\\
&\le \Pi \cdot d_{g_H}(p,q) + 2 \cdot l \cdot (d_{g_H}(p,q) + d_{g_H}(p,p_K)+d_{g_H}(q,q_K))\\
&\le (\Pi+ 4 \cdot l) \cdot d_{g_H}(p,q).
\end{align*}
Now, for the doubled twisted cone condition, we take a point $z$ on $\Gamma_{p_K,p}$ and observe that, along $\Gamma_{p,q}$, this point is closer to $p$ than to $q$. Hence, from \eqref{vk1}, $l_{min}(\Gamma_{p,q}(z)) = l_{g_H}(\Gamma_{p_K,p}(z)) \le 2 \cdot  \Pi/ \tau \cdot \delta_{\bp}(z)$,
and similarly for $z$ on $\Gamma_{q_K,q}$. Finally, for $z$ on $\gamma_{p_K,q_K}$ we first consider the subcase $z \in B_{1/8}(p_K) \cap \gamma_{p_K,q_K}$. The Lipschitz continuity of $\delta_{\bp}$ (assuming the Lipschitz constant to be $1$) shows that $\delta_{\bp}(p_K)/2 \le \delta_{\bp}(z)$, whence
\begin{align*}
l_{min}(\Gamma_{p,q}(z))&= l_{min}(\gamma_{{p_K,q_K}}(z)) + l_{g_H}(\Gamma_{p_K,p})\\
&\le  l/s \cdot \delta_{\bp}(z) + 2 \cdot  \Pi/ \tau \cdot \delta_{\bp}(p_K)\\
&\le (l/s + 4 \cdot  \Pi/ \tau) \cdot \delta_{\bp}(z).
\end{align*}
We treat the case $z \in  B_{1/8}(q_K) \cap \gamma_{{p_K,q_K}}$ in the same way. For the remaining case where
$z \in H \setminus (B_{1/8}(p_K) \cup B_{1/8}(q_K)) \cap \gamma_{p_K,q_K}$ we obtain $1/8 \le  l_{min}(\gamma_{{p_K,q_K}}(z)) \le  l/s \cdot \delta_{\bp}(z)$ from \eqref{pre}. Thus $l_{g_H}(\Gamma_{p_K,p}) \le  \Pi \cdot d_{g_H}(p,p_K) \le \Pi/4 \le 2 \cdot l/s \cdot \Pi \cdot  \delta_{\bp}(z)$,
and consequently $l_{min}(\Gamma_{p,q}(z)) = l_{min}(\gamma_{p_K,q_K}(z)) + l_{g_H}(\Gamma_{p_K,p}) \le  (1+2 \cdot \Pi) \cdot l/s \cdot \delta_{\bp}(z)$. Summarizing, $H \setminus \Sigma$ is $c$-\si-uniform for $c:= 4 \cdot(l + l/s + \Pi + \Pi/ \tau + \Pi \cdot l/s)$.\\

As in Step 3a above, we note that any two points $p,q \in H$ can be joined by some $c$-\si-uniform curve, supported in $H \cup \{p,q \}$, we find that $H$ is rectifiably connected. Since $H$ is compact the estimate $l_{g_H}(\Gamma(z)) \le  \Pi \cdot d_{g_M}(p,z)$ in \eqref{vk001} also shows that $H$ a finite intrinsic diameter: $diam_{g_H}H < \infty$.\qed

\noindent\textbf{Step 4 (Proof of the corollaries)}\\
To prove Corollary~\ref{fini} we use \eqref{vk001} which says that the intrinsic length  $l_{g_H}(\Gamma(z))$  of the pipeline $\Gamma(z)$ is upper bounded by the extrinsic distance of its endpoints $\Pi \cdot d_{g_{\R^{n+1}}}(p,z)$. (As in Step 3A, we note that we can choose one fixed $\varpi$ for all $H\in{\cal{H}}^{\R}_n$. Then $\Pi=\Pi(t_{n,\varpi},n)$ only depends on $n$.)  Since $l_{g_H}(\Gamma(z))$  is an upper bound for the intrinsic distance and the constructed $c$-\si-uniform curves  are compositions of two such pipelines we get a constant $c^{\R}_n \in (0,1)$ depending only on the dimension $n$ and such that $c^{\R}_n \cdot  d_{g_H}(p,q) \le  d_{g_{\R^{n+1}}}(p,q)$. The second inequality, $d_{g_{\R^{n+1}}}(p,q) \le d_{g_H}(p,q)$, is trivial.\\

For Corollary \ref{insu} we first claim that for any given $\ve >0$ and any two points $p,q \in F \setminus \Sigma_F$, there is a rectifiable curve $\gamma_{p,q}: [a,b] \ra F \setminus \Sigma_F$ with $\gamma_{p,q}(a)=p$ and $\gamma_{p,q}(b)=q$ which satisfies
\[
l_{g_F}(\gamma_{p,q})\le (c+\ve ) \cdot d_{g_F}(p,q)\quad\mm{and}\quad l_{min}(\gamma_{p,q}(z)) \le (c+\ve ) \cdot \delta_{\bp}(z)
\]
for any $z \in \gamma_{p,q}$. This follows from the fact that any compact subset of $F \setminus \Sigma_F$ admits arbitrarily fine $C^3$-approximations by suitable compact subsets of $k \cdot H$ for $k$ large enough. As in Step 2 we can then infer the existence of $\gamma_{p,q}$ from a corresponding curve in $k \cdot H$. Namely, we scale $H$ by a sufficiently large constant so that not only the $\D$-images of $p$ and $q$ can be identified with points in $F$ but also the $c$-\si-uniform curve that joins them in $H \setminus \Sigma$. This $\D$-preimage is a $(c+\ve )$-\si-uniform curve once we scaled $H$ appropriately, depending on the chosen $p$, $q$ and $\ve >0$. Finally, we send $\ve \ra 0$ and apply suitable BV-compactness results, namely Helly's selection principle \cite[Theorem 4 in Section 4.5]{SG}, to get a sequence $\gamma_n$ of $c+1/n$-\si-uniform curves subconverging to  a limit curve which is $c$-\si-uniform.
\qed
%
%
%
%
%
%%%%%%%%%%%%%%%%%%%%%%%%%%%%%%%%%
%%%%%%%%%%%%%%%%%%%%%%%%%%%%%%%%%
%%%%%%%%%%%%%%%%%%%%%%%%%%%%%%%%%
\setcounter{section}{3}
\renewcommand{\thesubsection}{\thesection}
\subsection{Hyperbolic Unfoldings}\label{hu}
%%%%%%%%%%%%%%%%%%%%%%%%%%%%%%%%%
%%%%%%%%%%%%%%%%%%%%%%%%%%%%%%%%%
%%%%%%%%%%%%%%%%%%%%%%%%%%%%%%%%%
On $H\setminus\Sigma$ we study the \si-metric $d_{\bp}$ as well as the quasi-conformal metric $k_{H\setminus\Sigma}$. Both metrics may be regarded as generalizations of the quasi-conformal metric $k_D$ on uniform domains $D\subset \R^n$. We discuss the resulting hyperbolic properties and determine their ideal Gromov boundary. A general reference on Gromov hyperbolic spaces is~\cite[Chapter III.H]{BH}.
%
%
%
%%%%%%%%%%%%%%%%%%%%%%%%%%%%%%%%%
\subsubsection{Quasi-Hyperbolic and \si-Metrics}\label{discrete}
%%%%%%%%%%%%%%%%%%%%%%%%%%%%%%%%%
Let $X$ be a non-complete, locally compact and complete, rectifiably connected metric space.

\begin{definition}\label{qh} \,
We define the \textbf{quasi-hyperbolic metric} $k_X$ for any two points $x,y\in X$ by
\[
k_X(x,y) := \inf\Bigl\{\int_\gamma 1/dist(\cdot,\p X) \, \, \Big| \, \gamma \subset  X  \mbox{ rectifiable curve  joining }  x \mbox{ and } y \Bigr\}.
\]
\end{definition}

From Proposition \ref{coh} we know that $H \setminus \Sigma$ is rectifiably connected. Thus we get the quasi-hyperbolic metric $k_{H \setminus \Sigma}$ on $X=H \setminus \Sigma$ with $\p X = \overline{X} \setminus X = \Sigma$. This metric uses only the intrinsic distance induced by $g_H$. With an \si-transform at hand we can build a new metric which encapsulates both information from the intrinsic metric $g_H$ as well as from its second fundamental form $A_H$:

\begin{definition}\label{sg} \,
For a given \si-transform $\bp$ we define on $H\setminus\Sigma$ the \textbf{\si-metric} $d_{\bp}$ by
\[
d_{\bp}(x,y) := \inf\Bigl\{\int_\gamma \bp \, \, \Big| \, \gamma   \subset  H
\setminus \Sigma\mbox{ rectifiable curve joining }  x \mbox{ and } y  \Bigr \}.
\]
 The length of a curve $\gamma \subset H\setminus\Sigma$ measured relative $d_{\bp}$ is denoted by $l_{\bp}(\gamma)$.
\end{definition}

\begin{remark}
In the definitions of $k_{H \setminus \Sigma}$ and $d_{\bp}$ we may actually allow curves in $H$. This follows from the path connectedness of $H \setminus \Sigma$ and the inequality $\bp(x) \ge L  / dist_{g_H}(x,\Sigma)$, since the expression $1/dist_{g_H}(x,\Sigma)$ considered along any curve reaching $\Sigma$ becomes non-integrable. In particular, such curves do not affect the infimum.
\end{remark}

Recall that a \textbf{geodesic curve}, or \textbf{geodesic} for short, is an isometric embedding $\gamma:[0,l]\subset\R\ra X$. A metric space is \textbf{geodesic} if any two points can be joined by a geodesic.

\begin{proposition}[Basic properties of $\mathbf{d_{\bp}}$ and $\mathbf{k_{H\setminus\Sigma}}$]\label{sgp}\hfill\newline\vspace{-16pt}
\begin{enumerate}
  \item Both $(H \setminus \Sigma, d_{\bp})$ and $(H \setminus \Sigma, k_{H \setminus \Sigma})$ are complete geodesic metric spaces.
  \item For any two $x$, $y \in H \setminus \Sigma$ we have
\[
d_{\bp}(x,y)\ge L \cdot k_{H\setminus \Sigma} (x,y)\ge  L/2 \cdot  log\Bigl(\bigl(1+\frac {d_{g_H}(x,y)}{dist_{g_H}(x, \Sigma_H)}\bigr)\cdot\bigl(1+\frac {d_{g_H}(x,y)}{dist_{g_H}(y, \Sigma_H)}\bigr)\Bigr),
\]
where $L$ denotes the Lipschitz constant for the \si-distance $\delta_{\bp}$.
  \item $(H \setminus \Sigma, d_{\bp})$ has \textbf{bounded geometry}*.
  \item For a flat norm converging sequence of minimizing hypersurfaces $H_i \ra H$, the \si-metrics $d_{\bp_{H_i}}$ converge compactly on smooth domains via $\D$-maps to $d_{\bp_H}$.
\end{enumerate}
\end{proposition}

\begin{remark}\label{dbg} \,
1.\ *The condition for bounded geometry is this. For global Lipschitz constant $l \ge 1$ and radius $\rho>0$ there exists around any point $p\in H\setminus\Sigma$ an $l$-bi-Lipschitz chart $\phi_p:B_\rho(p) \ra U_p$ between the ball $B_\rho(p)$ in $(H\setminus\Sigma, d_{\bp})$ to some open set $U_p \subset(\R^n,g_{\R^{n}})$. We shall always assume that $0 \in U_p$ and $\phi_{p}(p) = 0$. In cases where we need to specify these parameters we say that $(H\setminus\Sigma, d_{\bp})$  has  \emph{$(\rho,l)$-bounded geometry}.\\
2.\ The bounded geometry condition usually fails for $k_{H \setminus \Sigma}$ since there is no upper bound for $|A|(x) \cdot dist_{g_H}(x, \Sigma)$, $x \in H \setminus \Sigma$. In particular, there is no positive radius $\varrho>0$ that puts uniform constraints on the geometry of \emph{all} balls of radius $\varrho$ in $(H \setminus \Sigma, k_{H \setminus \Sigma})$.
\end{remark}

{\bf Proof of \ref{sgp}} \,
We choose $x,y \in H \setminus \Sigma$ and a smooth curve $\gamma :[0,1]\rightarrow H \setminus \Sigma$ with $\gamma (0)=x$ and $\gamma
(1)=y$. From $dist_{g_H}(x, \Sigma)+ d_{g_H}(x,\gamma (t)) \ge dist_{g_H}(\gamma (t), \Sigma)$ and $|\nabla \gamma(t)| \ge |\nabla d_{g_H}(x,\gamma(t))|$ we get
\[
\int\limits_\gamma \bp  \ge  L \cdot \int_\gamma \frac {1}{dist_{g_H}(\cdot,\Sigma)}\ge L \cdot \int\limits_0 ^1\frac {ds(t)}{dist_{g_H}(x, \Sigma)+
d_{g_H}(x,\gamma (t))} \ge L \cdot  log\left(1+\frac {d_{g_H}(x,y)}{dist_{g_H}(x, \Sigma_H)}\right),
\]
since $\bp(x) \ge L / dist_{g_H}(x,\Sigma)$. Similarly, we get the inequalities for $y$ instead of $x$. Adding these yields
\[
2\int\limits_\gamma \bp \ge 2L \cdot \int_\gamma \frac {1}{dist_{g_H}(\cdot,\Sigma)} \ge L \cdot \left( log\bigl(1+\frac
{d_{g_H}(x,y)}{dist_{g_H}(x, \Sigma_H)}\bigr) + log\bigl(1+\frac {d_{g_H}(x,y)}{dist_{g_H}(y, \Sigma_H)}\bigr) \right).
\]
This holds for all curves $\gamma$ connecting $x$ and $y$. Hence, we may pass to the infima $k_{H\setminus \Sigma}(x,y)$ and $d_{\bp}(x,y)$. This inequality also shows that both metrics are complete on $H \setminus \Sigma$. That both spaces are geodesic now easily follows from the Lipschitz continuity of $1/dist_{g_H}(x, \Sigma_H)$ and $\delta_{\bp}$ using Helly's selection principle \cite[Theorem 4 in Section 4.5]{SG} (compare also \cite[Chapter 10]{BHK}). Indeed, the spaces are complete, and for any two points with connecting curve $\gamma$ of upper bounded length we can find a compact subset so that $\gamma$ stays inside.\\

Now we turn to the proof of properties (iii) and (iv). For property (iii), we denote for $s\geq1$ the exponential map of $s\cdot M$ at $p$ by $\exp_p[s \cdot M]: (T_pM,g_{T_pM}) \ra s \cdot M$. Since $H$ is compact, standard results for the exponential map and the regularity theory of $H$ imply the following (cf.\ Step 2 in the proof of Proposition \ref{besi} for details): For each positive $\eta >0$ there is some $\Lambda(\eta) \gg 1$ so that for any $p \in H \setminus \Sigma$, the map $\exp_p[\Lambda(\eta) \cdot \bp(p) \cdot H]$ is a local diffeomorphism from $B_{10^3/L}(0)$ onto its image in $H$ with
\[
|\exp_p[\Lambda(\eta) \cdot \bp(p) \cdot H ]^*(\Lambda(\eta) \cdot \bp(p)^2 \cdot g_H) - g_{T_p H}|_{C^5(B_{10/L}(p))} \le \eta,
\]
where the ball $B_{10/L}(p)$ is taken with respect to the intrinsic metric of $\Lambda(\eta) \cdot \bp(p) \cdot H$.
Thus the radius is $10/(L \cdot \Lambda(\eta)) \ll 1/L$ relative $\bp(p)\cdot H$.\\

For small $\eta >0$, the exponential map is therefore $l$-bi-Lipschitz on $B_{10/(L \cdot \Lambda(\eta))}(p) \subset \bp(p) \cdot H$ for some $l\equiv l(\eta) \ra 1$ as $\eta \ra 0$. Note that, unlike the $C^5$-norm, the Lipschitz constant is invariant under scalings.\\

Property (iv) follows from the naturality of $\bp$ and the fact that the convergence upgrades to compact $C^k$-convergence for any $k \ge 0$.
\qed

\begin{corollary}\label{lreu} \, There are constants $\rho_n, l_n >0$ depending only of the dimension $n$ so that $(H \setminus \Sigma, d_{\bp})$ has $(\rho_n, l_n )$-bounded geometry, for any $H \in {\cal{H}}^{\R}_n$.
\end{corollary}

\noindent\textbf{Proof} \,
This still follows from the argument for Proposition \ref{sgp} (iii) using additionally, Proposition \ref{besi} for $H \in {\cal{H}}^{\R}_n$.  \ref{besi} gives the asserted uniform control for all points in $Q_H = H \setminus \Sigma$ and depending only on the dimension.
\qed

In Appendix~\ref{smsk} we explain how a Whitney type smoothing process can be applied to any \si-transform $\bp$. It generates a smooth $\bp^*>0$ on $H \setminus \Sigma$ such that for $\delta_{\bp^*}=1/\bp^*$ we have constants $c_i >0$, $i=1$, $2$, $3$, with
\begin{equation}\label{wsm}
c_1 \cdot \delta_{\bp}(x) \le \delta_{\bp^*}(x)  \le c_2 \cdot \delta_{\bp}(x)\quad\mm{and}\quad |\p^\beta \delta_{\bp^*}  / \p x^\beta |(x) \le c_3(\beta) \cdot \delta_{\bp}^{1-|\beta|}(x)
\end{equation}
for any $x \in H \setminus \Sigma$. (Here, the derivatives with multi-index $\beta$ are taken with respect to normal coordinates around $x \in H \setminus \Sigma$.) We set
\[
d_{\bp^*}(x,y) := \inf \Bigl  \{\int_\gamma \bp^* \, \, \Big| \, \gamma   \subset  H
\setminus \Sigma\mbox{ rectifiable curve joining }  x \mbox{ and } y  \Bigr \}
\]
This is the distance metric for the smooth Riemannian manifold $(H \setminus \Sigma, (\bp^*)^2 \cdot g_H)$.

\begin{corollary}\label{smskc} \,
The Whitney smoothing $(H \setminus \Sigma, d_{\bp^*})$ viewed as $(H \setminus \Sigma, (\bp^*)^2 \cdot g_H)$ is a complete Riemannian manifold with bounded geometry. It is quasi-isometric to $(H \setminus \Sigma, d_{\bp})$.\\

 Moreover, there are constants $\rho^*_n, l^*_n >0$ depending only of the dimension $n$ so that $(H \setminus \Sigma, d_{\bp^*})$ has $(\rho^*_n, l^*_n )$-bounded geometry, for any $H \in {\cal{H}}^{\R}_n$.
\end{corollary}

\noindent\textbf{Proof} \,
The first part is an immediate consequence of Proposition \ref{sgp} (iii) and \eqref{wsm}. The assertion for $H \in {\cal{H}}^{\R}_n$ follows when we additionally use \ref{lreu} and note that the Whitney smoothing constants $c_i$ of \ref{smsk} (ii) now merely depend on the dimension.
\qed
%
%
%
%%%%%%%%%%%%%%%%%%%%%%%%%%%%%%%%%
\subsubsection{Gromov hyperbolicity of $H\setminus\Sigma$}\label{gh}
%%%%%%%%%%%%%%%%%%%%%%%%%%%%%%%%%
If $\Sigma\n$, the uniformity of $H\setminus\Sigma$ (for $\Sigma\n$) immediately implies that
$(H \setminus \Sigma, k_{H \setminus \Sigma})$ is Gromov hyperbolic by the results of~\cite{GO} and~\cite[Theorem 3.6]{BHK}. Similarly, we will show that the \si-uniformity of $H\Sigma$ (where possibly $\Sigma\v$) makes the \si-metric $d_{\bp}$ on $H \setminus \Sigma$ again Gromov hyperbolic. We begin by recalling the

\begin{definition}[Gromov hyperbolicity]\label{grh} \,
A geodesic metric space is \textbf{Gromov hyperbolic}, or more precisely, $\mathbf{\delta}$\textbf{-hyperbolic}, if all its geodesic triangles are $\mathbf{\delta}$\textbf{-thin} for some $\delta >0$. This means that each point in an edge of a  geodesic triangle lies within $\delta$-distance of one of the other two edges. A complete Gromov hyperbolic space $X$ is called \textbf{visual}, or more precisely $\mathbf{\beta}$\textbf{-roughly starlike}, for $\beta>0$,  with respect to a base point $p\in X$, if for any $x \in X$ there is a geodesic ray starting at $p$ whose distance to $x$ is at most $\beta$.
\end{definition}

The concept of Gromov hyperbolic spaces is designed to study the asymptotic behavior near infinity. It embraces a broad range of spaces including objects like trees.

\begin{example}[Uniformity and hyperbolicity]\label{grex} \,
We describe some examples E1-E3 and counterexamples C1-C3 of spaces with hyperbolic properties related to our situation.
\begin{itemize}
  \item \textbf{E1} \, Compact Riemannian manifolds are always Gromov hyperbolic. We just choose $\delta = diameter$ and find that the manifold is ${\delta}$-hyperbolic. Similarly, we observe that for any compact Riemannian manifold $M$ the product $M \times \R$ is again Gromov hyperbolic.
  \item \textbf{E2} \, Consider an $H \in {\cal{H}}^c_n$ whose singulary set $\Sigma$ is a closed connected manifold. Assume that there is a neighborhood $U$ isometric to $\Sigma\times B_C$, where $B_C=B_1(0) \cap C$ for a complete connected area minimizing cone $C$ which is singular only in its tip at $0$. Then $(U, k_{H \setminus \Sigma})$ is a warped product:
  $B_C \setminus \{0\}$ is stretched to (one half of) an infinite cylinder, and the length of the $\Sigma$-fibers grows exponentially while we approach $0 \in B_C$. From this, the hyperbolicity of $k_{H \setminus \Sigma}$ can be check fairly directly.
  \item \textbf{E3} \, Euclidean uniform domains $D \subset \R^n$, and more generally, uniform spaces $X$ equipped with their quasi-hyperbolic metric $k_X$ on $X$ are Gromov hyperbolic, see~\cite{GO} and~\cite{BHK}. This example is universal in the sense that the uniformization theory of Bonk, Heinonen and Koskela~\cite{BHK} establishes a bijective conformal correspondence between the quasi-isometry classes of  \emph{proper geodesic and roughly  starlike  Gromov  hyperbolic} spaces on one hand side and the quasi-similarity classes of \emph{bounded locally compact uniform} spaces on the other. (A metric space is called \emph{proper} if all closed balls are compact.)
\end{itemize}
There are also well-known spaces where we easily find large and non-thin triangles so that these spaces are not Gromov hyperbolic.
\begin{itemize}
  \item \textbf{C1}\, Asymptotically flat spaces such as the Euclidean space are not Gromov hyperbolic, nor are products of non-compact complete Gromov hyperbolic spaces.
  \item \textbf{C2}\, Manifolds with sectional curvature $\equiv -1$ may not be Gromov hyperbolic, e.g., $\Z^2$-coverings of Riemann surfaces of genus $\ge 2$.
  \item \textbf{C3}\, For a compact manifold $(M^n,g_M)$, the product space $\R^{\ge 0} \times M^n$ equipped with the warped product metric $g_{\R} + (1+ a \cdot r)^2 \cdot g_M$, $a \ge 0$, is Gromov hyperbolic if and only if $a=0$.
\end{itemize}
Thus, although hyperbolicity is associated with fast growth of lengths and volumina, further spreading of Gromov hyperbolic spaces towards infinity can destroy their hyperbolicity.
\end{example}

Now we will see that $H \setminus \Sigma$ admits natural hyperbolic geometries with varying additional properties. For completeness we state the following result for the quasi-hyperbolic metric $k_{H\setminus\Sigma}$ which is due to~\cite{GO} and~\cite[Theorem 3.6]{BHK}.

\begin{proposition}\label{c22} \,
For any $H \in \cal{H}$ with $\Sigma\n$ we have
\begin{itemize}
  \item $(H \setminus \Sigma, k_{H \setminus \Sigma})$ is a complete Gromov hyperbolic space.
  \item If $H$ is compact, then $(H \setminus \Sigma, k_{H \setminus \Sigma})$ is roughly starlike.
\end{itemize}
\end{proposition}

However, unlike the quasi-hyperbolic metric on uniform Euclidean domains the Gromov hyperbolic space $(H \setminus \Sigma, k_{H \setminus \Sigma})$ need not to be of bounded geometry. This is a serious drawback in view of analytical arguments. For instance, we cannot expect uniform Harnack inequalities for elliptic problems. The geometric source for this non-boundedness are quickly sharpening \emph{wrinkles} in $H \setminus \Sigma$ corresponding to singular rays in the tangent cones of $H$, cf. Appendix \ref{flat-norm-approx}. Indeed, $\bp$ grows much faster than $1/dist_{g_H}(\cdot, \Sigma)$ when approaching $\Sigma$ along such wrinkles.\\

At any rate, it is natural to think of $d_{\bp}$ as a version of $k_{H \setminus \Sigma}$ which spreads further out near the boundary $\Sigma$. However, the metric space $(H \setminus \Sigma, k_{H \setminus \Sigma})$ (in general very roughly) resembles the direct product of E1 in Paragraph \ref{grex} near the boundary. Thus, in view of counterexample C3, a conformal deformation of $(H \setminus \Sigma, k_{H \setminus \Sigma})$ using $\bp \ge (L \cdot dist_{g_H}(\cdot,\Sigma))^{-1}$ with $\bp \gg (L \cdot dist_{g_H}(\cdot, \Sigma))^{-1}$ in wrinkled regions could potentially destroy the hyperbolicity of $k_{H \setminus\Sigma}$. It is precisely the sharper \si-uniformity of $H \setminus \Sigma$ which counterbalances this adverse spreading effect.

\begin{proposition}\label{c23} \,
Let $\bp$ be an \si-transform and $H\in{\cal{H}}_n$. Then both $(H \setminus \Sigma, d_{\bp})$ and its smoothing $(H \setminus \Sigma, d_{\bp^*})$ are complete Gromov hyperbolic spaces of bounded geometry. More precisely, if $H$ is $a$-\si-uniform, then
\begin{itemize}
  \item The \si-metric  $d_{\bp}$  is $\delta$-hyperbolic with $\delta=\delta(a, L_{\bp})$. For $L_{\bp}=1$, we have
\[
    \delta(a, 1)=4 \cdot a^2 \cdot \log \big(1+  c(a) \cdot ( 2 \cdot c(a) + 3)\big).
  \]
  \item The Whitney type smoothed \si-metric $d_{\bp^*}$ is $\Delta$-hyperbolic with $\Delta=\Delta(a,L_{\bp},H)$.
  \item For $H\in{\cal{H}}^\R_n$  we even have: $\delta(L_{\bp},n)$  and $\Delta(L_{\bp},n)$ independent of $H$.
\end{itemize}
\end{proposition}

\noindent\textbf{Proof} \,
We showed in Proposition~\ref{sgp} that  $(H \setminus \Sigma, d_{\bp})$ and $(H \setminus \Sigma, d_{\bp^*})$ are complete and have bounded geometry. The key ingredient for the hyperbolicity of $(H \setminus \Sigma, d_{\bp})$ is the \si-uniformity of $H \setminus \Sigma$. Nevertheless, our proof is modelled on the strategy for proving the hyperbolicity of $(D,k_D)$ for uniform domains $D \subset \R^n$, see Gehring and Osgood~\cite{GO} and Bonk, Heinonen, Koskela~\cite[Chapter 2-3]{BHK}.\\

For the remainder of this proof, we assume that $H \setminus \Sigma$ is $a$-\si-uniform for some fixed $a \ge 1$. For ease of notation we assume that  $L_{\bp}=1$, i.e., $|\delta_{\bp}(x)-\delta_{\bp}(y)| \le d_{g_H}(x,y)$ for  $x, y \in H \setminus \Sigma$. The case $L_{\bp}\neq 1$ follows similarly. We subdivide the proof into three lemmas.

\begin{lemma}[Relations between $\mathbf{d_{g_H}}$ and $\mathbf{d_{\bp}}$]\label{ue} \,
For any two $x$, $y \in H \setminus \Sigma$ we have
\begin{equation}\label{ri}
d_{\bp}(x,y) \le 4 \cdot a^2 \cdot \log \big(1+d_{g_H}(x,y) \cdot\max \{\bp(x), \bp(y)\}\big),
\end{equation}
and in particular
\begin{equation}\label{risq}
d_{\bp}(x,y) \le   4 \cdot a^2 \cdot \sqrt{d_{g_H}(x,y) \cdot\max \{\bp(x), \bp(y)\}.}
\end{equation}
Conversely, we have
\begin{equation}\label{ril}
\log\big(1+l_{g_H}(\gamma(x,y)) \cdot \max \{\bp(x), \bp(y)\}\big)  \le d_{\bp}(x,y),
\end{equation}
where $l_{g_H}(\gamma(x,y))$ is the length of the $d_{\bp}$-geodesic curve $\gamma(x,y)$ in $H\setminus\Sigma$ measured with respect to the intrinsic distance induced by $g_H$. In particular, we obtain
\begin{equation}\label{riw}
\log\big(1+d_{g_H}(x,y) \cdot \max \{\bp(x), \bp(y)\}\big)\le d_{\bp}(x,y)
\end{equation}
and
\begin{equation}\label{riww}
\big|\log\delta_{\bp}(x) - \log\delta_{\bp}(y)\big| \le d_{\bp}(x,y).
\end{equation}
\end{lemma}

\noindent\textbf{Proof} \,
Let $\gamma \subset H \setminus \Sigma$ be an $a$-\si-uniform curve joining $x$ and $y \in H \setminus \Sigma$ which is of length $\lambda= l_{g_H}(\gamma)$. Choose the midpoint $z \in \gamma$, i.e., $\gamma = \gamma_1 \cup \gamma_2$ with $\{z\} = \gamma_1 \cap \gamma_2$ and $l_{g_H}(\gamma_1)=l_{g_H}(\gamma_2)$ for the subcurves $\gamma_i \subset \gamma$, $x \in \gamma_1$, $y \in \gamma_2$. We claim the following  inequalities:
\begin{equation}\label{ri1}
l_{\bp}(\gamma_1) \le 2 \cdot a \cdot \log(1+ \lambda\cdot \bp(x))\quad\mm{and}\quad  l_{\bp}(\gamma_2) \le 2 \cdot a \cdot \log(1+ \lambda \cdot \bp(y)),
\end{equation}
where $l_{\bp}$ denotes the length with respect to $d_{\bp}$. We first use the Lipschitz estimate \eqref{li} for $\delta_{\bp}$ to prove that if $l_{g_H}(\gamma_1) = \lambda/2 < \delta_{\bp}(x)$, then
\begin{equation}\label{ri3}
l_{\bp}(\gamma_1) = \int_{\gamma_1} 1/\delta_{\bp} \le \int_0^{\lambda/2}1/(\delta_{\bp}(x) - s) \, ds.
\end{equation}
Indeed, after parameterizing $\gamma_1$ by arc length, \eqref{li} shows that $\delta_{\bp}(\gamma(s)) \ge \delta_{\bp}(x) - s$ when we leave $x$ at time $0$, whence \eqref{ri3}. %We use this estimate when either we actually know that $l_{g_H}(\gamma_1) = \lambda/2 < \delta_{\bp}(x)$ (in case A below) or we apply it to a subcurve of $\gamma_1$ that satisfies this condition (in case B below).
We now distinguish the following two cases:
\[
\mm{\textbf{A.}} \;\; l_{g_H}(\gamma_1) \le \frac{a}{a+1} \cdot \delta_{\bp}(x)\quad\mm{and}\quad  \mm{\textbf{B.}}\;\; l_{g_H}(\gamma_1) > \frac{a}{a+1} \cdot \delta_{\bp}(x).
\]
\noindent\textbf{A.} \, From \eqref{ri3} we infer
\begin{align*}
l_{\bp}(\gamma_1) &\le \int_0^{\lambda/2}\frac{1}{\delta_{\bp}(x) - s}\, ds = \log\left(\delta_{\bp}(x)/(\delta_{\bp}(x)-\lambda/2)\right)\\
&= \log\left(1/(1-\lambda/2\delta_{\bp}(x))\right) \le 2a \cdot \log(1 + \lambda/\delta_{\bp}(x)).
\end{align*}
Here we used the elementary inequality $\log(1/(1-x)) \le 2k \cdot \log(1+2x)$ for any $k \ge 1$, $x \in [0,k/(k+1)]$.\\

\noindent\textbf{B.} \, The $a$-\si-uniformity shows that $a/l_{g_H}(\gamma_1[0,t]) \ge 1/
\delta_{\bp}(\gamma_1(t))$ for $t \le \lambda/2$. We combine this inequality with \eqref{ri3} for the subcurve of $\gamma_1$ from $x$ to the point where the length attains the value $\frac{a}{a+1} \cdot \delta_{\bp}(x)$. Since $a \ge 1$, we get
\begin{align*}
l_{\bp}(\gamma_1) &\le \int_0^{\frac{a}{a+1} \cdot
\delta_{\bp}(x)}1/(\delta_{\bp}(x) - s) \; ds+  a \cdot \int_{\frac{a}{a+1} \cdot \delta_{\bp}(x)}^{\lambda/2}1/s \; ds\\
&=\log\Big(\frac{1}{1-\frac{a}{a+1}}\Big) + a \cdot \log\Big(\frac{(a+1) \cdot \lambda}{2 \cdot a \cdot \delta_{\bp}(x)}\Big) \le \log(1+a)+a \cdot \log\Big(\frac{\lambda}{\delta_{\bp}(x)}\Big)\\
&\le a \cdot \log 2 + a \cdot \log (1+\lambda/\delta_{\bp}(x)) \le 2 \cdot a \cdot \log (1+\lambda/\delta_{\bp}(x)),
\end{align*}
where we used that $\lambda > 2 \cdot \frac{a}{a+1} \cdot \delta_{\bp}(x) \ge \delta_{\bp}(x)$ and applied the elementary inequality
\begin{equation}\label{ri2}
\log(1+ k \cdot x) \le k \cdot \log(1+   x) \mm{ for } k \ge 1,\, x \ge 0.
\end{equation}
Thus for both cases A and B the first inequality of \eqref{ri1} holds; the second is established similarly. Moreover, the inequality \eqref{ri} follows from $l_{g_H}(\gamma) \le a \cdot d_{g_H}(x,y)$, and another application of \eqref{ri2}. As we remarked above this implies \eqref{risq}.\\

For~\eqref{ril}, we choose $x$, $y \in H \setminus \Sigma$ and consider a rectifiable curve $\gamma$ in $(H \setminus \Sigma, d_{\bp})$ with $\gamma (0)=x$ and $\gamma (1)=y$. By Lipschitz continuity we have $\delta_{\bp}(y)\le\delta_{\bp}(x) + d_{g_H}(x,y)$ which gives $\delta_{\bp}(\gamma (t)) \le \delta_{\bp}(x)+ d_{g_H}(x,\gamma (t))$ and $\delta_{\bp}(\gamma (t)) \le \delta_{\bp}(x)+ l_{g_H}(x,\gamma (t))$,
where $l_{g_H}(x,\gamma (t))$ is the length of the subcurve $\gamma([0,t])$ measured in $(H \setminus \Sigma, g_H)$. From this, we note the inequalities
\[
\frac {\delta_{\bp}(\gamma (t))}{\delta_{\bp}(x)} \le \frac {\delta_{\bp}(x)+d_{g_H}(x,\gamma (t))}{\delta_{\bp}(x)} \le \frac {\delta_{\bp}(x)+l_{g_H}(x,\gamma (t))}{\delta_{\bp}(x)}
\]
In the case where $\gamma$ is a geodesic curve in $(H \setminus \Sigma, d_{\bp})$ we find
\[
\log\left(\frac {\delta_{\bp}(x)+l_{g_H}(\gamma(x,y))}{\delta_{\bp}(x)}\right)
= \int\limits_0 ^1\frac {ds(t)}{\delta_{\bp}(x)+ d_{g_H}(x,\gamma (t))} \le \int\limits_0 ^1\frac {ds(t)}{\delta_{\bp}(\gamma (t))} \le d_{\bp}(x,y).
\]
Exchanging the r\^oles of $x$ and $y$ we get the same inequalities for $y$ instead of $x$. From both sets of inequalities, as well as $d_{g_H}(x,y) \le l_{g_H}(\gamma(x,y))$, we deduce \eqref{ril}, \eqref{riw} and \eqref{riww} using the elementary inequality $\log(1+x) \le \sqrt{x}$, $x \ge 0$.
\qed

\begin{lemma}[\si-uniformity of geodesic curves]\label{hul} \,
Every geodesic curve $\gamma$ in $(H \setminus \Sigma, d_{\bp})$ is a \textbf{c-\si-uniform} curve in $(H \setminus \Sigma, g_H)$ for some $c\ge a$ independent of $\gamma$.
\end{lemma}

\noindent\textbf{Proof} \,
Let $\gamma$ in $(H \setminus \Sigma, d_{\bp})$ be a geodesic curve between $x$ and $y\in H \setminus \Sigma$. We check the two defining conditions for a c-\si-uniform curve.\\

\emph{Twisted double cone conditions} \,
We set $D:= \max_{z \in \gamma} \delta_{\bp}(z)$. Then we can find unique integers $N(x)$, $N(y) \ge 0$ with
\[
D/2^{N(x)+1} < \delta_{\bp}(x) \le D/2^{N(x)}\quad\mm{and}\quad D/2^{N(y)+1} < \delta_{\bp}(y) \le D/2^{N(y)}.
\]
Next we subdivide $\gamma$ into subcurves. Choose points $x_0,\ldots,x_{N(x)}$ and $y_0,\ldots,y_{N(y)}\in\gamma$ such that by starting from $x$, $x_i$ is the first point on $\gamma$ where $\delta_{\bp}(x_i) = D/2^i$. Since $\delta_{\bp}$ is continuous, $x_i$ obvioiusly exists. Similarly we define the points $y_j$ starting from $y$. This defines geodesic curves $\gamma_x(i)$ between $x_i$ and $x_{i+1}$, and $\gamma_y(j)$ between $y_j$ and $y_{j+1}$, as well as a curve $\gamma_0$ between $x_0$ and $y_0$. Since $\bp = 1/\delta_{\bp}$ we have $\bp(\gamma_x(i)) \subset  [2^i/D,\infty)$, $\bp(\gamma_y(j)) \subset  [2^j/D,\infty)$ and $\bp(\gamma_0) \subset  [1/D,\infty)$. Consequently, since $\gamma$ is a geodesic, \eqref{risq} of Lemma~\ref{ue} shows that
\begin{itemize}
  \item $l_{g_H}(\gamma_0)/D \le  l_{\bp}(\gamma_0) \le 4a^2 \cdot \sqrt{2 \cdot l_{g_H}(\gamma_0)/D}$.
  \item $l_{g_H}(\gamma_x(i))/(D/2^{i}) \le  l_{\bp}(\gamma_x(i)) \le 4a^2 \cdot \sqrt{l_{g_H}(\gamma_x(i))/(D/2^{i+1})}$.
\end{itemize}
Similarly, we get the analogous estimates for $\gamma_y(j)$. From these inequalities we deduce
\begin{equation}\label{rr}
l_{g_H}(\gamma_0)/D \le  32a^4,\quad l_{g_H}(\gamma_x(i))/(D/2^{i}) \le 32a^4,\quad l_{g_H}(\gamma_y(j))/(D/2^{j}) \le 32a^4
\end{equation}
which in turn implies $l_{\bp}(\gamma_0) \le 32a^4$, $l_{\bp}(\gamma_x(i)) \le 32a^4$, and $l_{\bp}(\gamma_y(j))  \le  32a^4$. We  use this to estimate $\bp$ on $\gamma_x(i)$ away from its endpoints. Towards this end, we recall that in the starting and end points $x_i$ and $x_{i+1}$ of $\gamma_x(i)$ we
have $\delta_{\bp}(x_i)=D/2^i$ and $\delta_{\bp}(x_{i+1})=D/2^{i+1}$. The same holds for $\gamma_y(i)$. Next let $z \in \gamma_x(i)$. From \eqref{riw} in Lemma~\ref{ue} we deduce
\[
\big|\log(D/2^i) - \log(\delta_{\bp}(z)) \big|\le d_{\bp}(x_i,z) \le l_{\bp}(\gamma_x(i)) \le 32a^4
\]
and therefore  $\exp(-32a^4) \cdot D / 2^i \le \delta_{\bp}(z)$. Using \eqref{rr}, we finally arrive at   \[
l_{min}(\gamma(z))  \le \sum_{k \ge i} 32a^4 \cdot D \cdot 2^{-k}  \le 64a ^4 \cdot D / 2^i \le b(a) \cdot \delta_{\bp}(z)
\]
for $b(a):= 64a^4 \cdot  \exp(32 a^4)$.\\

\emph{Quasi-geodesics} \,
On $\gamma$ we choose two points $\tilde x$ and $\tilde y$ so that for the subcurves $\gamma_{\tilde x}$ from $x$ to $\tilde x$ and $\gamma_{\tilde y}$ from $y$ to $\tilde y$ we have $l_{g_H}(\gamma_{\tilde x})=l_{g_H}(\gamma_{\tilde y})=d_{g_H}(x,y)/2$.  Then each of these curves reaches at most the midpoint of $\gamma$ whence
\begin{itemize}
  \item the length of the curve between $\tilde x$ and $\tilde y$ is $l_{g_H}(\gamma(\tilde x,\tilde y))=l_{g_H}(\gamma)-d_{g_H}(x,y)$.
  \item $l_{g_H}(\gamma_{\tilde x}) \le  b(a) \cdot  \delta_{\bp}(\tilde x)$ and $l_{g_H}(\gamma_{\tilde y}) \le  b(a) \cdot \delta_{\bp}(\tilde y)$ with $b(a)$ as defined above.
  \item $d_{g_H}(\tilde x,\tilde y) \le 2 \cdot d_{g_H}(x,y)$ by the triangle inequality.
\end{itemize}
Now \eqref{riw} of Lemma \eqref{ue} gives
\begin{align*}
&\log\big(1+ (l_{g_H}(\gamma) - d_{g_H}(x,y)) \cdot  \max\{\bp(\tilde x),\bp(\tilde y)\}\big)\le d_{\bp}(\tilde x,\tilde y)\\
&\le 4a^2 \cdot \log\big(1+2 \cdot d_{g_H}(x,y)\cdot\max\{\bp(\tilde x),\bp(\tilde y)\}\big) \le 4a^2 \cdot \log(1+ 4b(a))
\end{align*}
Thus for $b^*\equiv b^*(a):= \exp\big(4a^2 \cdot \log(1+ 4b)\big)-1$ we have
\begin{equation}\label{ror}
l_{g_H}(\gamma) \le d_{g_H}(x,y) + b^*\cdot \min\{\delta_{\bp}(\tilde x),\delta_{\bp}(\tilde y)\}.
\end{equation}
Next we distinguish the cases
\[
\mathbf{A.} \;\min\{\delta_{\bp}(\tilde x),\delta_{\bp}(\tilde y)\} \le (4a^2+1) \cdot d_{g_H}(x,y)\quad\mm{and}\quad\mathbf{B.} \mm{ otherwise.}
\]
Combining case \textbf{A} with \eqref{ror} yields $l_{g_H}(\gamma) \le\big(1 + b^*\cdot (4a^2+1)\big) \cdot d_{g_H}(x,y)$. In case \textbf{B} we get
\[
(4a^2 + 1) \cdot d_{g_H}(x,y) \le \min\{\delta_{\bp}(\tilde x),\delta_{\bp}(\tilde y)\}\le\min\{\delta_{\bp}(x),\delta_{\bp}(y)\} +  d_{g_H}(x,y)
\]
using again the Lipschitz condition on $\delta_{\bp}$. This also implies $d_{g_H}(x,y)/\min\{\delta_{\bp}(x),\delta_{\bp}(y)\} \le 1/4a^2$. Consequently, if $a \ge1$,
\[
d_{\bp}(x,y) \le  4a^2 \cdot \log \big(1+d_{g_H}(x,y) \cdot \max \{\bp(x), \bp(y)\}  \big) \le 1
\]
since $x \ge \log(1+x)$ for $x > 0$. Now \eqref{ri} and \eqref{ril} from Lemma~\ref{ue} show that
\[
\log (1+ l_{g_H}(\gamma) \cdot \max \{\bp(x), \bp(y)\}) ) \le d_{\bp}(x,y) \le 1.
\]
From this and $\log(2) >1/2$, $x \le 2 \cdot \log(1+x)$ for $x \in [0, 1]$, and $x \ge \log(1+x)$ for $x > 0$,
we see
\begin{align*}
l_{g_H}(\gamma) \cdot \max \{\bp(x), \bp(y)\})&\le 2\log\big(1+ l_{g_H}(\gamma) \cdot \max \{\bp(x), \bp(y)\}\big) \le 2d_{\bp}(x,y)\\
&\le 8a^2 \cdot \log \big(1+d_{g_H}(x,y) \cdot \max \{\bp(x), \bp(y)\} \big)\\
&\le 8a^2 \cdot  d_{g_H}(x,y) \cdot \max \{\bp(x), \bp(y)\}
\end{align*}
using again \eqref{ri} and \eqref{ril} from Lemma~\ref{ue}. Summarizing, we get $l_{g_H}(\gamma) \le 8a^2 \cdot  d_{g_H}(x,y)$ in case B. Then $c(a):= 1 + b^*\cdot (4a^2+1) +  8a^2$ does the job for both A and B. We conclude that $\gamma$ is a $c(a)$-\si-uniform curve.
\qed

\begin{lemma}[Geodesic triangles are thin]\label{tg} \,
There is some $\delta(a) >0$ so that every geodesic triangle in $(H\setminus\Sigma, d_{\bp})$ is $\delta$-thin.
\end{lemma}

\noindent\textbf{Proof} \,
Let $x$, $y$ and $z \in H \setminus \Sigma$ be the vertices of a geodesic triangle in $(H \setminus \Sigma, d_{\bp})$, and let $[x,y]$, $[y,z]$ and $[x,z]$ be the three geodesic edges. These are $c$-\si-uniform for some $c\ge a \ge 1$ by Lemma \ref{hul}. We claim that there is some $\delta\equiv\delta(a) >0$ such that
\[
d_{\bp}(p,[y,z] \cup [x,z]) \le \delta
\]
for any $p \in [x,y]$. For the sake of concreteness, assume $l_{g_H}([x,p]) \le l_{g_H}([y,p])$. \si-uniformity implies
\begin{equation}\label{tt}
l_{g_H}([x,p]) \le c \cdot \delta_{\bp}(p)\quad\mm{ and }\quad l_{g_H}([x,y]) \le  c \cdot d_{g_H}(x,y).
\end{equation}
Now we distinguish the cases
\[
\mathbf{A.} \;\; c \cdot l_{g_H}([x,z])  <  l_{g_H}([x,p])\quad\mm{ and }\quad\mathbf{B.} \;\; c \cdot l_{g_H}([x,z])  \ge  l_{g_H}([x,p]).
\]
In case \textbf{A}, we know from the assumption  $l_{g_H}([x,p]) \le l_{g_H}([y,p])$ that $2 \cdot l_{g_H}([x,p]) \le l_{g_H}([x,y])$ and, thus, in case A, we have $2 \cdot l_{g_H}([x,z]) \le l_{g_H}([x,y])$, since $c\ge 1$. The triangle inequality then shows that $2 \cdot l_{g_H}([y,z]) \ge l_{g_H}([x,y]) \ge  2 \cdot  l_{g_H}([x,p])$. Therefore, there is some $q \in [y,z]$ with $l_{g_H}([q,z])  =  (2c)^{-1} \cdot l_{g_H}([x,p]) \le l_{g_H}([y,q])$ and
\[
d_{g_H}(p,q) \le l_{g_H}([x,p])+l_{g_H}([x,z])+l_{g_H}([q,z]) \le (1 + c^{-1} + c^{-1}/2) \cdot l_{g_H}([x,p]).
\]
The $c$-\si-uniformity shows that $l_{g_H}([x,p])/2 = c \cdot l_{g_H}([q,z]) \le c^2 \cdot \delta_{\bp}(q)$. Furthermore, \eqref{tt} gives $l_{g_H}([x,p])/2 = c/2 \cdot \delta_{\bp}(p)$, whence
\[
d_{g_H}(p,q) \le 2c^2 \cdot (1 + c^{-1} + c^{-1}/2) \cdot \min \{\delta_{\bp}(p), \delta_{\bp}(q)\}.
\]
In case \textbf{B}, there is some $q \in [x,z]$ with $l_{g_H}([q,x])  =  (2c)^{-1} \cdot l_{g_H}([x,p]) \le l_{g_H}([z,q])$ and
\[
d_{g_H}(p,q) \le l_{g_H}([x,p])+l_{g_H}([x,q]) \le (1 + c^{-1}/2) \cdot l_{g_H}([x,p]).
\]
The $c$-\si-uniformity shows that $l_{g_H}([x,p])/2 = c \cdot l_{g_H}([q,x]) \le c^2 \cdot \delta_{\bp}(q)$. Using $l_{g_H}([x,p])/2 = c/2 \cdot \delta_{\bp}(p)$ again yields
\[
d_{g_H}(p,q) \le 2 \cdot c^2 \cdot (1 +  c^{-1}/2) \cdot \min \{\delta_{\bp}(p), \delta_{\bp}(q)\}.
\]
Summarizing, we get $d_{g_H}(p,q) \le 2 \cdot c^2 \cdot (1 + c^{-1} + c^{-1}/2) \cdot \min \{\delta_{\bp}(p), \delta_{\bp}(q)\}$ in both cases. This, the inequality
\[
d_{\bp}(p,[y,z] \cup [x,z]) \le  d_{\bp}(p,q) \le 4 \cdot a^2 \cdot \log\big(1+d_{g_H}(p,q) /\min \{\delta_{\bp}(p), \delta_{\bp}(q)\}\big),
\]
and \eqref{ri} in Lemma~\ref{ue} finally imply
\[
d_{\bp}(p,[y,z] \cup [x,z]) \le 4a^2 \cdot \log \big(1+  c(a) \cdot ( 2c(a) + 3)\big)=:\delta(a).
\]
Put differently, every geodesic triangle in $(H \setminus \Sigma, d_{\bp})$ is $\delta(a)$-thin, hence $(H \setminus \Sigma, d_{\bp})$ is $\delta(a)$-hyperbolic.
\qed

To finish the proof we note that Gromov hyperbolicity is a quasi-isometric invariant~\cite[Chapter III.H (1.9)]{BH}. Hence, the Whitney smoothing $(H \setminus \Sigma, d_{\bp^*})$ of $(H \setminus \Sigma, d_{\bp})$ is Gromov hyperbolic by Corollary~\ref{smskc}.  [BH], Ch.III.H (1.7)-(1.9) show that   $(H \setminus \Sigma, d_{\bp^*})$ is $\Delta(L_{\bp},H)$-hyperbolic. Finally, for $H \in {\cal{H}}^{\R}_n$ we already know from \ref{intsuni} that $H \setminus \Sigma$ is $c_n$-\si-uniform and using  \ref{smsk}(ii) we see $\Delta(L_{\bp},H)=\Delta(L_{\bp},n)$,.\qed

\begin{remark}[\si-Metrics on Regular Spaces]\label{sire} \,
1.\ The classical quasi-hyperbolic metric $k_{H \setminus \Sigma}$ has no meaning when $H$ is regular. The natural convention, namely $dist_{g_H}(\cdot,\Sigma)\equiv\infty$, would lead to $k_{H \setminus \Sigma} \equiv 0$ so that $diam(H,k_{H \setminus \Sigma})=0$. Hence $(H,k_{H \setminus \Sigma})$ would degenerate to a one-point space. On the other hand, the metric $d_{\bp}$ is well-defined and non-trivial even if $H$ is a regular hypersurface. For $\Sigma_H \v$, the metric space $(H,d_{\bp})$ is homeomorphic to $(H,g_H)$, unless $H$ is totally geodesic where again $diam(H, d_{\bp})=0$. On the other hand, if $H$ is compact then $diam(H, d_{\bp})= \infty$ if and only if $H$ is singular. In this way, we can think of the diameter $diam(H, d_{\bp})$ as a measure for the \emph{relative curvature} of $H$ inside $M$.\\
2.\ If a sequence of compact area minimizing hypersurfaces $H_i$ converges in flat norm to a compact minimizer $H_\infty$, then, by (S4) (naturality) $\bp_{H_i}$ converges in $C^\alpha$-norm to $\bp_{H_\infty}$, $\alpha \in (0,1)$, near any given regular point of $H_\infty$. Hence, via $\D$-maps, we have compact convergence
\[
(H_i\setminus \Sigma_{H_i},d_{\bp_{H_i}}) \ra (H_\infty \setminus \Sigma_{H_\infty}, d_{\bp_{H_\infty}}) \quad\mm{ as }i \ra \infty.
\]
3.\ The various \si-metrics in a given converging sequence $(H_i\setminus \Sigma_{H_i},d_{\bp_{H_i}})$ are not only individually Gromov hyperbolic, but they are all $\delta$-\emph{hyperbolic for the same} $\delta >0$. To see this, we first note that the hyperbolicity constant of $(H_i \setminus \Sigma_{H_i},d_{\bp_{H_i}})$ only depends on the \si-uniformity parameter. But we know from Proposition~\ref{intsuni}~(iv) that there is a common constant $c_n$, for all $H \in {\cal{H}}^{\R}_n$, for which $H \setminus \Sigma$ is $c_n$-\si-uniform. The \si-uniformity constant for $H_i\in{\cal{H}}_n$ can then be estimated in terms of this $c_n$ and the constant for regions away of $\Sigma_i$, which in turn are controlled by $H_\infty$. This argument extends to the Whitney smoothings $(H_i \setminus \Sigma_{H_i},d_{\bp^*_{H_i}})$ from an additional use of the estimates in~\cite[Chapter III.H (1.7)-(1.9)]{BH}.\\
4. We note in passing that similar phenomena happen for degenerating families of smooth Riemann surfaces of genus $\geq2$ equipped with their hyperbolic metric. We have smooth convergence of these metrics to the limit metric in smooth regions, whereas the family will develop infinite complete ends where the limit surface has singular points.
\end{remark}
%
%
%
%%%%%%%%%%%%%%%%%%%%%%%%%%%%%%%%%
\subsubsection{$\Sigma \subset H$ as a Gromov Boundary}\label{grch}
%%%%%%%%%%%%%%%%%%%%%%%%%%%%%%%%%
We use the hyperbolicity of $(H \setminus \Sigma,d_{\bp})$ and $(H \setminus \Sigma, k_{H \setminus \Sigma})$ to describe the singular set $\Sigma\n$ as an ideal boundary for a particular compactification of these spaces.\\

\noindent\textbf{Basic concepts.} \,
Let $X$ be a complete Gromov hyperbolic space. A \textbf{geodesic ray} is an isometric embedding $\gamma:[0,\infty)\ra X$. A \textbf{generalized geodesic ray} $\gamma:I\ra X$ is either a geodesic curve or a geodesic ray. In the former case where $I=[0,l]$ we extend $\gamma$ to a ray by defining $\gamma(t)=\gamma(l)$ for $t\in[l,\infty)$. Two geodesic rays are equivalent if they have finite Hausdorff distance. The equivalence class of a ray $\gamma$ will be denoted by $\gamma(\infty)$.

\begin{definition}[Gromov boundary] \,
The set $\p_G X$ of equivalence classes of geodesic rays is called the \textbf{Gromov boundary} of $X$.
\end{definition}

Using the extension of a geodesic curve to a geodesic ray we can identify $\overline{X}_G = X \cup \p_G X$ with $\{\gamma(\infty)|\,\gamma\mm{ is a generalized ray}\}$. Moreover, given $q\in X$ any equivalence classe $\gamma(\infty)$ may be represented by a geodesic ray starting at $q$~\cite[Lemma III.H.3.1]{BH}. We define a topology on $\overline{X}_G$ as follows. We say that a sequence $x_n \in \overline{X}$ converges to $x \in \overline{X}$ if there exist generalized rays $c_n$ with $c_n(0) = q$ and $c_n(\infty) = x_n$ subconverging (on compact sets) to a generalized ray $c$ with $c(0) = q$ and $c(\infty) = x$. Then $\p_G X$ is closed, $\overline{X}_G$ is compact, and the canonical map $X\hookrightarrow\overline{X}_G$ is a homeomorphism onto its image,~\cite[Proposition III.H.3.7]{BH}. Furthermore, $\overline{X}_G$ is metrisable~\cite[Chapter III.H.3]{BH}. It is called the \textbf{Gromov compactification of $X$}.\\

\noindent\textbf{Identification of $\p_GX$.} \, For the flat model of a uniform domain $D \subset \R^n$, the Gromov boundary of the complete space $X=(D,k_D)$ is well-understood: There is a canonical bijection between $\p_GX$ and $\p D$ which assigns to each geodesic ray in $X$ its end point in $\p D$~\cite[Theorem 3.6]{BHK}.\\

The counterparts for the three complete spaces $X_{\bp}:= (H \setminus \Sigma,d_{\bp})$, its Whitney smoothing $X_{\bp^*}:= (H \setminus \Sigma,d_{\bp^*})$ and $X_{1/dist}:= (H \setminus \Sigma, k_{H \setminus \Sigma})$ read as follows.

\begin{theorem}\label{xgb} \,
For singular $H \in {\cal{H}}^c_n$, the identity map on $H \setminus \Sigma$ extends to a homeomorphism between $H$ and the Gromov compactifications of $X_{\bp}$, $X_{\bp^*}$ and $X_{1/dist}$:
\[
H \cong(\overline{X}_{\bp})_G \cong (\overline{X}_{\bp^*})_G \cong (\overline{X}_{1/dist})_G,
\]
where $\cong$ means homeomorphic. In particular, we have $\Sigma \cong\p_GX_{\bp} \cong \p_GX_{\bp^*} \cong \p_GX_{1/dist}$.\\

For singular $H \in {\cal{H}}^{\R}_n$ the identity map on $H \setminus \Sigma$ extends to a homeomorphism between the one-point compactification $\widehat{H}$ of $H$ and the Gromov compactifications of $X_{\bp}$, $X_{\bp^*}$ and $X_{1/dist}$:
\[
\widehat{H} \cong(\overline{X}_{\bp})_G \cong (\overline{X}_{\bp^*})_G \cong(\overline{X}_{1/dist})_G.
\]
In particular, we have $\widehat{\Sigma} \cong\p_GX_{\bp} \cong \p_GX_{\bp^*} \cong \p_GX_{1/dist}$.
\end{theorem}

\noindent\textbf{Proof} \,
For $X_{1/dist}$, that is, the uniform space $H \setminus \Sigma$ equipped with its quasi-hyperbolic metric $k_{H \setminus \Sigma}$, the result follows from the general theory of uniform spaces~\cite[Theorem 3.6 and Proposition 3.12] {BHK} and the definition of the topology for the Gromov compactification. The case of the \si-metrics can be treated in a quite similar way.\\

In essence the idea is this: For $H \in {\cal{H}}^c_n$ we fix a base point $p\in X$ and assign to (equivalence classes of) geodesic rays in $X_{\bp}$ starting at $p$ their end point which actually lies in $\Sigma\subset H$. For $H \in {\cal{H}}^{\R}_n$ we also have \si-uniform curves of \emph{infinite} length with respect to $(H \setminus \Sigma, g_H)$. These account for the point at infinity of the one-point compactification. Therefore, we start with the case $H \in {\cal{H}}^c_n$ before extending the argument to hypersurfaces in ${\cal{H}}^{\R}_n$.\\

\noindent\textbf{Case A: $H \in {\cal{H}}^c_n$.} \,
For $H \in {\cal{H}}^c_n$ we define a \textbf{canonical bijection} $\Psi_{\Sigma}:\p_GX_{\bp} \ra \Sigma$.\\

Towards this end let $\gamma:[0,L) \ra H \setminus \Sigma$, $L \in (0,\infty]$, be a proper geodesic ray in $X_{\bp}$ starting from $p\in H\setminus \Sigma$ which relative to $(H \setminus \Sigma, g_H)$ is parameterized by arc-length and has length $L$. From Lemma \ref{hul}, $\gamma$ is a $c$-\si-uniform curve for some $c(H) >0$. Thus, since $H$ is compact and $diam\,X_{\bp} < \infty$, the quasi-geodesic condition for $\gamma$ shows that $L < \infty$. We claim that for $t < L$, $t \ra L$, there exists a point $x\in\Sigma$ such that $\gamma(t) \ra x$. Indeed, since $[0,L]$ is the maximal interval of definition, there must be a sequence $t_i \in (0,L)$, $t_i \ra L$ as $i \ra \infty$, so that $\gamma(t_i) \ra x$ for some $x \in \Sigma$. Moreover, the quasi-geodesic condition on $\gamma$ implies that $\gamma(s_i) \ra x$ for any other sequence $s_i \in (0,L)$ with $s_i \ra L$.\\

Next consider two such geodesic rays  $\gamma[1]$ and  $\gamma[2]$ with end points $x[k] \in \Sigma$, and which have finite Hausdorff distance in $X_{\bp}$, that is, they define
the same point in $\p_GX_{\bp}$. Then we find sequences $t_i[k] \in\big(0,L(\gamma[k])\big)$ with $t_i[k]\ra L(\gamma[k])$, $k=1,2$, so that $d_{\bp}(t_i[1],t_i[2]) \le c = const< \infty$. Further, we note that $\bp(t_i[k]) \ra \infty$ as $i \ra \infty$. From \eqref{riw}, that is, $\log\big(1+d_{g_H}(x,y) \cdot \max \{\bp(x), \bp(y)\}\big)\le d_{\bp}(x,y)$ we have $d_{g_H}(t_i[1],t_i[2]) \ra 0$ as $i \ra \infty$, whence
$x[1]=x[2]$. Thus every representative of a point in $\p_GX_{\bp}$ has the same endpoint in $\Sigma$. This yields a well-defined map $\Psi_{\Sigma}$ from $\p_GX_{\bp}$ to $\Sigma$. We claim that $\Psi_{\Sigma}$ is \textbf{bijective}.\\

\noindent\textbf{Surjectivity of $\Psi_{\Sigma}$.} \,
Let $x \in \Sigma$. We choose a sequence $x_i \in H \setminus \Sigma$ with $x_i\ra x$ as $i \ra \infty$, and a sequence of geodesic curves $\gamma_i$ from $p$ to $x_i$. Then, using the Arzel\`{a}-Ascoli theorem, we get a compactly converging subsequence of the $\gamma_i$ with limiting geodesic $\gamma$. From the previous argument we see that $\gamma$ links $p$ with some $y \in \Sigma$. The quasi-geodesic condition on $\gamma_i$ shows then that $y=x$.\\

\noindent\textbf{Injectivity of $\Psi_{\Sigma}$.} \,
For geodesic rays $\gamma[1]$ and  $\gamma[2]$ with end points $x[1]=x[2] \in \Sigma$ we choose two sequences $t_i[k] \in\big(0,L(\gamma[k])\big)$, $i=0,1,\ldots$, with
\begin{equation}\label{di}
t_i[k] \ra L(\gamma[k]),\, k=1,2\quad\mm{and}\quad l_i:=L(\gamma[1])-t_i[1] = L(\gamma[2])-t_i[2].
\end{equation}
From the \si-uniformity we infer that for large $i \gg 1$,
\[
\bp\big(\gamma[k](t_i[k])\big) \le  2c /\big(L(\gamma[k])-t_i[k]\big).
\]
In turn, the triangle inequality shows that
\[
d_{H}\big(\gamma[1](t_i[1]),\gamma[2](t_i[2])\big) \le L(\gamma[1])-t_i[1] + L(\gamma[2])-t_i[2]=2l_i.
\]
Since
\[
d_{\bp}\big(\gamma[1](t_i[1]),\gamma[2](t_i[2])\big)\le 4a^2 \cdot \sqrt{2l_i \cdot  2c /l_i}  \le   8a^2 \cdot \sqrt{c},
\]
using \eqref{risq} from Lemma \ref{ue} and \eqref{di} implies that $d_{\bp}\big(\gamma[1](t_i[1]),\gamma[2](t_i[2])\big)$ remains bounded as $i\ra \infty$. From this we infer that $\gamma[1]$ and  $\gamma[2]$ have finite Hausdorff distance in $X_{\bp}$ and thus determine the same point in $\p_GX_{\bp}$.\\

To conclude Case A, we first note that $\Psi_\Sigma$ is continuous by a proof along the lines of the surjectivity of $\Psi_\Sigma$. Since $\Psi_\Sigma$ is bijective, $\p_GX_{\bp}$ and $\Sigma$ are compact and $\overline{X_{\bp}}_G$ is metrizable, $\Psi_{\Sigma}$ must be a homeomorphism. Summarizing, we see that the map $\Phi_H: \overline{X_{\bp}}_G  \ra H$ which is defined by $\Phi_H|_ {H \setminus \Sigma}=id_ {H \setminus \Sigma}$ and $\Phi_H|_{\p_GX_{\bp}}=\Psi_{\Sigma}$ yields a \textbf{homeomorphism} extending the identity on $H \setminus \Sigma$.\\

Finally, for $X_{\bp^*}$ we use that $(H \setminus \Sigma, d_{\bp^*})$ and $(H \setminus \Sigma, d_{\bp})$ are quasi-isometric. In particular, their Gromov compactifications and boundaries are homeomorphic, cf.~\cite[Theorem III.H.3.9]{BH}.\\

\noindent\textbf{Case B: $H \in {\cal{H}}^\R_n$.} \,
We take again a proper geodesic ray $\gamma:[0,L) \ra H \setminus \Sigma$, $L \in (0,\infty]$ in $X_{\bp}$ starting at $p$. Further, $\gamma$ is parameterized by arc-length and is of length $L$ relative to $(H \setminus \Sigma, g_H)$. This time we either have $L < \infty$ or $L=\infty$.\\

For $L < \infty$ we can argue as in Case A and get a homeomorphism $\Psi^*_{\Sigma}$ from $\p^*_GX_{\bp}$ to $\Sigma$, where $\p^*_G X \subset \p_G X$ denotes the set of equivalence classes of geodesic rays with finite length relative to $(H \setminus \Sigma, g_H)$. Now in a given equivalence class of geodesic rays, the representing curves all have either finite or infinite length. Indeed, each subcurve is again $c$-\si-uniform. From the twisted double \si-cone condition we see that $t\le c \cdot \delta_{\bp}(\gamma(t))$ for any $t >0$ and geodesic ray $\gamma$ parametrized by arc-length and of infinite length relative to $(H \setminus \Sigma, g_H)$. On the other hand, for a geodesic ray $\gamma^*$ with infinite length relative to $(H \setminus \Sigma, g_H)$, inequality \eqref{risq} from Lemma~\ref{ue} asserts that
\[
d_{\bp}(\gamma(t),\gamma^*(t)) \le 4a^2 \cdot\big(d_{g_H}(\gamma(t),\gamma^*(t))/\min  \{\delta_{\bp}(\gamma(t)), \delta_{\bp}(\gamma^*(t))\}\big)^{1/2}
\]
while the triangle inequality gives $d_{g_H}(\gamma(t),\gamma^*(t)) \le 2t$. Since $d_{\bp}(\gamma(t),\gamma^*(t)) \le8a^2$ for any $t >0$, the geodesic rays $\gamma$ and $\gamma^*$ are equivalent. Conversely, for a geodesic ray $\gamma^*$ which determines the same point in the Gromov boundary as $\gamma$, \eqref{riww} in Lemma \ref{ue} shows that
\[
\big|\log(\delta_{\bp}(\gamma(t))) - \log(\delta_{\bp}(\gamma^*(t))) \big| \le d_{\bp}(\gamma(t),\gamma^*(t)).
\]
Hence $\gamma^*$ has infinite length relative to $(H \setminus \Sigma, g_H)$. Consequently, there is precisely one point $z_\infty$ in $\p_GX_{\bp}$ corresponding to geodesic curves with infinite length relative to $(H \setminus \Sigma, g_H)$. Any of these geodesic rays leaves any bounded set in $H$, since otherwise it would approach some $z \in \Sigma$. But these points are reached by rays of finite
length. Thus they all approach $\infty_H$ and we may identify $z_\infty$ with $\infty_H$. In conclusion, we can extend the homeomorphism $\Psi^*_{\Sigma}$ from $\p^*_GX_{\bp}$ to $\Sigma$ to a homeomorphism $\Psi_{\widehat{\Sigma}}$  from $\p_GX_{\bp}$ to $\widehat{\Sigma}$. The remaining assertions follow as in Case A.
\qed
\appendix
%
%
%
%%%%%%%%%%%%%%%%%%%%%%%%%%%%%%%%%
%%%%%%%%%%%%%%%%%%%%%%%%%%%%%%%%%
%%%%%%%%%%%%%%%%%%%%%%%%%%%%%%%%%
\setcounter{section}{1}
\renewcommand{\thesubsection}{\thesection}
\subsection{Oriented Boundaries and Currents}\label{all}
%%%%%%%%%%%%%%%%%%%%%%%%%%%%%%%%%
%%%%%%%%%%%%%%%%%%%%%%%%%%%%%%%%%
%%%%%%%%%%%%%%%%%%%%%%%%%%%%%%%%%
In this appendix we gather some ideas, concepts and notations from geometric measure theory for the case of (almost) area minimizing hypersurfaces.\\

\noindent\textbf{I.\ Existence of area minimizers.} \,
Here, a convenient tool is the theory of oriented (minimal) boundaries, see for instance~\cite{AFP}, \cite{Gi} and \cite{MM}. In the language of geometric measure theory these correspond to locally normal currents of codimenson $1$.\\

Let $\Omega\subset \R^{n+1}$ be a bounded open set, and $f \in L^1(\Omega,\R)$. We define
\[
\int_\Omega|D f| := \sup \{\int_\Omega f \cdot \mbox{div} g  \, d\mu \;| \; g \in C_0^1(\Omega,\R^{n+1}), |g|_{C^0} \le 1 \}.
\]
We call $f$ a function of \textbf{bounded variation} or just \textbf{BV-function in $\Omega$} if $\int_\Omega|D f| < \infty$. The set of BV-functions on $\Omega$ is denoted by $BV(\Omega)$. The \textbf{BV-norm} of a BV-function $f$ is defined as $|f|_{BV(\Omega)} := |f|_{L^1(\Omega)} + \int_\Omega|D f|$. Finally, $f\in BV_{loc}(\Omega)$ if $f\in BV(\Omega_0)$ for any $\Omega_0\Subset\Omega$, i.e., $\Omega_0$ is bounded with $\overline\Omega_0\subset\Omega$.\\

If $f=\chi_E$ is the characteristic function of some Borel set $E \subset \R^{n+1}$ one refers to $\int_\Omega |D \chi_E|$ as the \textbf{perimeter} $P(E,\Omega)$ of $E$ in $\Omega$, since $P(E,\Omega)$ equals the $n$-dimensional Hausdorff-measure ${\cal{H}}_n(\p E \cap \Omega)$ if the boundary $\p E$ is smooth, cf.\ \cite[Example 1.4]{Gi}. A Borel set with locally finite perimeter, that is, $P(E,\Omega)<\infty$ for all open $\Omega\Subset\R^{n+1}$, is called a \textbf{Caccioppoli set}.\\

If $\Omega\Subset\R^{n+1}$ is open and $L$ is a Caccioppoli set, we can find an area minimizing hypersurface $E$ with $E\setminus\Omega=L\setminus\Omega$ by taking a perimeter minimizing sequence $\chi_{E_j}$ of Caccioppoli sets $E_j$ with $E_j \setminus \Omega \equiv L \setminus \Omega$. For $\p \Omega$ sufficiently smooth, for instance if $\Omega$ has Lipschitz regular boundary, the embedding $BV_{loc}(\Omega) \hookrightarrow L^1_{loc}(\Omega)$ is compact~\cite[3.23]{AFP}. Hence there is a subsequence $E_{j_k}$ converging in $L_{loc}^1$~\cite[Theorem 1.19]{Gi}. By lower semicontinuity of BV-norms \cite[Theorem 1.9]{Gi}, the limit $E$ is again a Caccioppoli set.\\

\noindent\textbf{II.\ Regularity theory for almost minimizers.} \,
While existence of (almost) minimizers is rather straight forward, regularity issues are very intricate. De Giorgi and others developed a partial regularity theory for minimal Caccioppoli sets which was actually extended to the more general case of \textbf{almost minimizers} by Tamanini \cite{T1}, \cite{T2}, Massari and Miranda \cite{MM}, Bombieri \cite{Bo} and Allard \cite{A}. The following result is taken from \cite[Theorem 1]{T1}.

\begin{definition}\label{defarm} \,
Let $\Omega \subset \R^{n+1}$ be open. The boundary $\p E$ of a Caccioppoli set $E \subset \R^{n+1}$ is called \textbf{almost minimizing in $\Omega$} if for some $K > 0$, $\alpha \in (0,1)$ and $R>0$, the inequality
\[
\psi(E, B_\rho(x)) := \int_{B_\rho(x)}|D \chi_{E}| - \inf \left\{\int_{B_\rho(x)}|D \chi_{F}|\, \Big| \, F
\Delta E \Subset B_\rho(x) \right\} \le K \cdot \rho^{n+2 \cdot \alpha}
\]
holds for any $x \in \Omega$ and $\rho \in (0,R)$ (here $F \Delta E := F \setminus E \cup  E \setminus F$).
\end{definition}

\textbf{Area minimizers} correspond to the case $\psi \equiv 0$. On the other hand, the hypersurfaces $S_C=\p B_1(0) \cap C$ obtained from an area minimizing cone $C \subset \R^{n+1}$ with tip at $0$ are almost minimizers.

\begin{proposition}\label{arm} \,
For an almost minimising boundary $\p E$ in $\Omega\subset\R^{n+1}$, $\p E \cap \Omega$ is a $C^{1, \alpha}$- hypersurface except for a singular set $\Sigma$ of Hausdorff codimension greater or equal than $8$.
\end{proposition}

Further improvements of regularity can be obtained from standard elliptic theory. For instance, the smooth locus of an \emph{area minimizing} hypersurface is analytic if the ambient manifold is analytic, see for instance~\cite[Chapter 5.7]{Mo}. This clearly holds in the case of Euclidean boundaries.\\

Using local coordinate charts these definitions and regularity results carry over to Riemannian manifolds without difficulties. Indeed, diffeomorphisms of the ambient space map preserve the condition of being an almost minimizer, for they locally preserve the estimate on $\psi(E, B_\rho(x))$ up to multiplication by the $n$th-power of the local maximum of the norm of their Jacobian. Thus in a Riemannian manifold, an almost minimizer is a hypersurface which via charts can be locally mapped to Euclidean almost minimizers in the sense of Definition~\ref{defarm}.\\

Proposition~\ref{arm} also implies that a sequence of almost minimizers $E_i$ converging to some limit $E_\infty$ will eventually become smooth near smooth limit points in $\p E_\infty$,~\cite[Theorem 1]{T1}. Further, $L^1_{loc}$-convergence implies $C^1$-convergence when the limit is known to be $C^{1,\alpha}$-smooth, see Allard's work~\cite{A} or Simon's lecture notes~\cite[Theorem 23.1]{Si2} for details.

\begin{corollary}\label{carm} \,
Let $E_i$, $i \ge 0$, be a sequence of almost minimizers for fixed $(K,\alpha)$ in some open bounded set $\Omega$.
\begin{enumerate}
  \item Assume that $E_i\to E_\infty$ in $L^1_{loc}$ with points $p_i \in \p E_i\to p_\infty \in \p E_\infty$. If $p_\infty$ is a smooth point in $\p E_\infty$, then so is, for sufficiently large $i$, the point $p_i\in\p E_i$.
  \item If the limit $E_\infty$ in (i) has a $C^{1, \alpha}$-boundary in $\Omega$, then $\p E_i$ converges to $\p E_\infty$ in $C^1$-topology.
\end{enumerate}
\end{corollary}

\begin{remark} \,
1.\ Note that $E_\infty$ also satisfies $\psi(E_\infty, B_\rho(x)) \le K \cdot \rho^{n+2 \cdot \alpha}$. This can be proved as in~\cite[Lemma 9.1]{Gi}.\\
2.\ Corollary~\ref{carm} carries over to Riemannian manifolds and asserts that a flat norm converging sequence of area minimizers will be locally $C^k$-converging around smooth points of the limit surface, cf.\ \cite[Lemma 11.4]{Gi} for details.
\end{remark}

A typical scenario for such convergence results are blow-ups at some $p\in \Sigma_H$ for a given $H \in {\cal{G}}$, that is, rescaling $H$ around $p$ by a sequence $\tau_m\to\infty$. Then there is a subconverging sequence $\tau_{m_k}\cdot H$ whose limit is an area minimizing cone. Formally this reads as follows cf.\ \cite[4.3.16]{F1}, \cite[Chapter 37.4]{Si1}, \cite[Theorem 1]{T1}.

\begin{proposition}\label{flat-norm-approx}
Let $H \in {\cal{G}}$ and $p \in \Sigma_H$. For every sequence $\tau_m \to +\infty$ of positive real numbers there exists a subsequence $\tau_{m_k}$, as well as an area minimizing cone $C_p \subset \R^{n+1}$, with $0 \in \sigma_{C_p}$, we call a \textbf{tangent cone}, such that
\begin{itemize}
  \item \textbf{flat norm convergence:} For any given open $U \subset \R^{n+1}$ with compact closure the \textbf{flat norm} $\db_U$ converges to zero: \, $\db_U (\tau_{m_k} \cdot H, C_p) \to 0.$
  \item  \textbf{$C^l$-norm convergence:} If, in addition, $\overline U \subset C_p \setminus \sigma_{C_p}$, then $\db_U$-convergence implies compact $C^l$-convergence, for any $l \ge 0$, expressible via  $\D$-maps, cf. Ch.\ref{nat},
\end{itemize}
$\sigma_{C_p}$ is our generic notation for singularities of cones.  $\db_U$ is defined in \eqref{fll}, roughly speaking, it measures the volume between $\tau_{m_k} \cdot H$ and $C_p$ in $U$.
\end{proposition}

We also note some well-known applications of the regularity theory for area minimizing hypersurfaces which, however, are hard to localize in the literature.

\begin{corollary}\label{t3} \,
Let $D \subset H$ be an open domain in an oriented minimal boundary $H$. Then the following statements are equivalent:
\begin{enumerate}
 \item All points in $D$ are regular.
 \item Near any point $p \in D$ the norm of the second fundamental form is bounded, i.e., $|A| \le c$ for some $c\equiv c(p) >0$.
 \item For any $p \in D$, the tangent cone is a hyperplane.
\end{enumerate}
\end{corollary}

Next we state a \emph{non-extinction} result for oriented minimal boundaries in $\R^{n+1}$ which is crucial for our compactness arguments. Roughly speaking it asserts that sequences of such minimizers cannot form approaching opposing sheets which annihilate in the limit.

\begin{lemma}\label{nex} \,
Let $H_i \subset \R^{n+1}$ be a sequence of oriented minimal boundaries with $0\in H_i$ and $|A_{H_i}| \le 1$ on $B_2(0) \subset H_i$. Then for any compactly converging subsequence $H_{i_k}$, the limit hypersurface $H_\infty$  is an oriented minimal boundary with $0 \in H_\infty$ and $|A_{H_\infty}| \le 1$ on $B_1(0) \subset H_\infty$. Furthermore, the $B_1(0) \subset H_{i_k}$ converge smoothly to $B_1(0) \subset H_\infty$ in the sense of $\D$-maps (see Definition~\ref{idmap}).
\end{lemma}

\noindent{\bf Proof} \,
We show that in $\R^{n+1}$ the ball $B_1(0) \subset H_i$ is not approached from $O_i:=H_i \setminus B_2(0)$ as $i \ra \infty$. This implies that there is a lower positive distance bound between $B_1(0)$ and the $O_i$, independent of $i$. Let us assume the contrary. As a standard consequence of DeGiorgi-Allard regularity theory~\cite[Theorem 24.2]{Si1} and the Harnack inequality~\cite[p.\ 73]{So} we could write a subset of $O_i$ as a smooth graph $G_i$ over $B_i:=B_1(0) \subset H_i$ which is arbitrarily close to $B_1(0)$ in $C^3$-norm for $i \gg 1$. Since the $H_i$ bound open sets $U_{H_i} \subset \R^{n+1}$ we may assume that $G_i$ and $B_i$ have opposite orientation. For $i \gg 1$ we consider $G_i \cup B_i$ and join $\p G_i$ and $\p B_i$ linearly through some hypersurface $F_i$. Then we add the bounded open set $V_i\subset \R^{n+1}$ with $\p V_i = G_i \cup B_i \cup F_i$ to $U_{H_i}$ and form the new open set $\tilde U_{H_i}=U_{H_i}\cup V_i\cup G_i\cup B_i$. For a sufficiently small bounded open set $\Omega$ containing $\overline V_i$ we have $\tilde U_{H_i} \cap \Omega^c=U_{H_i}$ while in $\Omega$, $\tilde U_{H_i}$ has smaller area than $H_i$ for $i\gg1$. (Note that the $B_i$ have uniformly bounded geometry from $|A_{H_i}|\big|_{B_2(0)} \le 1$.) This contradicts the area minimizing property of $H_i$ since $\tilde U_{H_i}$ is a compactly supported variation. Hence $O_i$ remains in a positively lower bounded distance of $B_i$ for all $i$. The remaining assertions follow from regularity theory.
\qed

We also note a weak Harnack type inequality for $|A|$.

\begin{lemma}\label{haa} \,
For any $\lambda\in (0,1]$ and $R_0 >0$ there is a constant $c\equiv c(\lambda, n, R_0) >0$ such that for any oriented minimal boundary $H \subset \R^{n+1}$ and $p \in H$ with $\sup \{|A|(x) \, | \, x \in B_{R_0}(p) \cap H\} \ge 1$, we have
\[
\sup \{|A|(x) \, | \, x \in B_{\lambda  \cdot R_0}(p) \cap H\} \ge c.
\]
\end{lemma}

\noindent\textbf{Proof} \,
Without loss of generality we may assume that $R_0 =1$ and $p=0 \in \R^{n+1}$. Assume that there is no such constant $c > 0$. Then there is some $\lambda \in (0,1)$ and a sequence of hypersurfaces $H_k$ with $0 \in H_k$ so that $\sup \{|A|(x) \, | \, x \in B_\lambda(0) \cap H_k\} \le 1/k$. Due to the minimality of these hypersurfaces there is a subsequence $H_{k_j}$ which on $\R^{n+1}$ converges compactly in flat norm to some limit hypersurface $H_\infty$. As in \ref{nex} we may assume from $\sup \{|A|(x) \, | \, x \in B_\lambda(0) \cap H_k\} \le 1/k$ that this is $C^k$-convergence in $B_\lambda(0)$ for some $k \ge 5$. Then the analytic minimizer $H_\infty$ is a hyperplane, since the limit of the $B_\lambda(0) \cap H_k$ in $H_\infty$ must be flat. Since $H \setminus \Sigma_H$ is connected by Proposition~\ref{coh}, the regularity theory promotes the flat convergence to $C^k$-convergence also outside $B_\lambda(0)$. But then $\sup \{|A|(x) \, | \, x \in B_R(0) \cap H_k\} \ra 0$ for any $R >0$, contradicting the assumption.
\qed

\noindent\textbf{III.\ \si-Structures on almost minimizers and Plateau solutions.} \,
The reasoning for our Theorems, as stated in Ch.\ref{int},  extends from the case of area minimizers in ${\cal{H}}$, we considered in Ch.\ref{siu} and  Ch.\ref{hu}, to the more general case of \textbf{almost minimizers} in ${\cal{G}}$. Most of the arguments carry over  to almost minimizers (and similarly to Plateau problems) unchanged. The few adjustments needed are discussed in following.\\

The Definition~\ref{lsk} of metric \si-transforms equally applies to  almost minimizers and the axioms (S1) - (S3) remain valid. The naturality condition (S4) still holds for
converging sequences of varifolds with commonly bounded generalized mean curvature. Then the same arguments as before, now based on Allard theory ~\cite[Ch.5]{Si1} apply. The case of blow-ups fits into this scenario. However, for the purposes of this paper we note that  the naturality axiom (S4) for blow-ups already follows from \ref{flat-norm-approx}.\\

The blow-up naturality of $\bp$, and not its broader variant on ${\cal{H}}$, is sufficient to establish Theorem~\ref{thm2} (\si-uniformity on $H\setminus\Sigma$). Namely, in Ch.\ref{siu} we only appeal to this form of naturality when we derive estimates from the limit which, even for almost minimizers, always belongs to ${\cal{H}}^{\R}_n$. In particular, the Bombieri-Giusti version of the localized isoperimetric inequality for oriented minimal boundaries~\cite[Theorem 2, p.\ 31]{BG} also applies to small balls on almost minimizers; their proof consists precisely in considering blow-up limits.\\

 Further, the hyperbolic unfolding Theorems~\ref{thm3} and~\ref{thm4} are based on \si-uniformity and not on the (almost) minimality. In Ch.\ref{hu}, where these Theorems are proved, we do not use that the \si-uniform spaces are (almost) minimizers. (The Remark \ref{sire} is not needed in the arguments but only describes some extensions and further context.) Thus, they also extend to almost minimizers and this even holds for the results on Whitney smoothings in Theorem~\ref{thm4} since the finer properties of \si-adapted covers in Proposition~\ref{besi} below again merely use that blow-up limits belong to ${\cal{H}}^{\R}_n$. Hence Theorem~\ref{thm3} holds for almost minimizers except that $d_{\bp}$ only commutes for blow-ups.\\

Next we turn to bounded area minimizers $H$ with boundary $\p H$, that is, $H$ solves the \textbf{Plateau problem} for the boundary $\p H$. For the sake of simplicity, we assume that $\p H$ is $C^2$-smooth. Due to boundary regularity results of~\cite{HS} this implies that $H$ is a $C^{1,\alpha}$-regular manifold near $\p H$. In particular, $\Sigma \cap \p H \v$.\\

 We replace $\delta_{\bp}(x)$ by $\mathbf{d}(x):=\min \{dist_{g_H}(x,\p H),\delta_{\bp}(x)\}$. Under the condition that the singular $H$ is a \emph{uniform space} one deduces a version of \si-uniformity as in Theorem~\ref{thm2} after replacing $\delta_{\bp}$ by $\mathbf{d}$. Then we construct hyperbolic unfoldings for $\mathbf{d}$ (and similarly for its Whitney smoothing $\mathbf{d}^*$) by merging the argument of Theorem~\ref{thm3} for $\p H$ and of Theorem~\ref{thm4} for $\Sigma$ for the distance function
\[
d_{\mathbf{d}}(x,y) := \inf \Bigl  \{\int_\gamma  1/\mathbf{d}(\cdot) \, \, \Big| \, \gamma   \subset  H \setminus \Sigma\mbox{ rectifiable curve joining }  x \mbox{ and } y  \Bigr \}.
\]
As in Chapter \ref{hu} we deduce that $d_{\mathbf{d}}$ and $d_{\mathbf{d}^*}$ define complete Gromov hyperbolic spaces with bounded geometry such that
\[
H \cong \overline{(H \setminus \Sigma,d_{\mathbf{d}})}_G  \cong \overline{(H \setminus \Sigma,d_{\mathbf{d}^*})}_G\quad\mm{and}\quad\Sigma \cup \p H \cong\p_G(H \setminus \Sigma,d_{\mathbf{d}}) \cong \p_G(H \setminus \Sigma,d_{\mathbf{d}^*}).
\]

\noindent\textbf{IV.\ Currents.} \,
So far we considered minimizers of codimension one which arose as boundaries of fairly general Borel sets. Next we want to enlarge the class of submanifolds of any codimensions via distributions. Let $U \subset \R^n$ be open and consider the $\mathcal{D}^m(U)$, the space of smooth $m$-forms compactly supported in $U$. This inherits a natural topology from the space of smooth functions compactly supported in $U$, and we consider its topological dual $\mathcal{D}_m(U)$, the space of $m$-\textbf{currents}. In particular, any submanifold $N$ of codimension $k$ defines an $n-k$-current $\llb N \rrb$ by integration.\\

We define the \textbf{boundary} $\p T  \in \mathcal{D}_{m}(U)$ of  $T \in \mathcal{D}_{m+1}(U)$ as the current  $\p T(\omega):= T(d \omega)$. The \textbf{support} $supp \: T $ of a current $T$ is the complement of the union of all open sets $W$ such that $T (\omega) = 0$ for  $\omega \in \mathcal{D}^n(U)$ with
$supp \; \omega \subset W$. For any open $W \subset U$ and $T  \in \mathcal{D}_{m}(U)$ we write $T \llcorner W$ for the current in $\mathcal{D}_{m}(W)$ we get from \textbf{restricting} $T$ to $\mathcal{D}_m(W)$. For any compactly supported current  $T  \in \mathcal{D}_{m}(U)$ we define its \textbf{push-forward} $f_\sharp T$ by $f_\sharp T(\omega):= T(f^*\omega)$, where $f^*\omega$ denotes the usual pull-back of the $m$-form $\omega$.\\

We let ${\bf{M}}_U (T) := \sup_{|\omega| \le 1, \mbox{{\tiny{supp}}} \omega \subset U}T(\omega)$ be the \textbf{mass} of the current $T$. To define the \textbf{flat (pseudo-)metric} on $\mathcal{D}_{m}(U)$ we consider open subsets $W \subset \overline{W} \subset U\subset \mathbb{R}^n$. Then
\begin{equation}\label{fll}
\db_W(C_1,C_2):=\mbox{inf}\{{\bf{M}}_W (S)+{\bf{M}}_W (R)\, |\, C_1 - C_2 = S + \p R, S \in \mathcal{D}_m(U), R \in \mathcal{D}_{m+1}(U)) \}.
\end{equation}
The family of these (pseudo-)metrics $\db_W$ generate the \textbf{flat norm topology} on $\mathcal{D}_m(U)$.\\

Finally, we define some important subclasses of currents. We call a current $T \in \mathcal{D}_m(U)$ \textbf{integer multiplicity rectifiable} or \textbf{rectifiable} for short, if for any $\ve >0$ and any compact set $K \subset U$ there exists a compactly supported $m$-dimensional polyhedral chain with $\Z$-coefficients of oriented simplices $P=P(K, T,\ve) \subset \R^k$ and a Lipschitz map $f:\R^k \ra \R^n$ such that $supp \, f_\sharp P \subset K$ and ${\bf{M}}_U (T-f_\sharp P) < \ve$. We denote by $\mathcal{R}_m(U) \subset \mathcal{D}_m(U)$ the space of integer multiplicity rectifiable currents and by $\mathcal{I}_m(U) \subset \mathcal{D}_m(U)$ the space of \textbf{integral currents}. Here, a current $T$ is integral if $T$ and $\p T$ are rectifiable currents with rectifiable boundary. There are compactness results for integral currents expressed in terms of the flat metric topology similar to compactness in $L^1$-topology for BV-functions. All these concepts and notions extend to compact manifolds via local charts. In particular, we get the following basic existence result, cf.\ \cite[4.2.17, 4.4.5 and 5.1.6]{F1} or \cite[Corollary 1 in 5.4.1]{GMS}, and \cite[Section 5.3]{F1} or \cite{F2} for the regularity assertions.

\begin{proposition}\label{hm} \,
For any $\alpha \in H_{n}(M^{n+1}, \Z)$ there is a mass minimizing integral current $X^n \in \alpha$ whose support is a smooth hypersurface outside a set $\Sigma_X \subset M^{n+1}$ of codimension greater or equal than $8$.
\end{proposition}

\noindent\textbf{V.\ Decomposition of rectifiable currents.} \,
To make contact with the theory of oriented boundaries we note the following decomposition theorem for rectifiable currents~\cite[4.5.17]{F1}, \cite[Theorem 7 in 4.3.1]{GMS} or \cite[Chapter 37]{Si1}.

\begin{proposition}\label{dic} \,
For any $R \in \mathcal{R}_n(\R^{n+1})$, with $\p R \v$, there exist measurable sets $A_i \subset \R^{n+1}$, $i \in \Z$ with $A_i \subset A_{i+1}$  such that for any bounded open $W \subset \R^{n+1}$, we have
\[
R= \sum_{i \in \Z} \p \llb A_i \rrb\quad\mm{and}\quad {\bf{M}}_W (R)= \sum_{i \in \Z} {\bf{M}}_W (\p \llb A_i \rrb).
\]
\end{proposition}

In the case of a locally mass minimizing current one may assume that the sets $A_i$ are open and the $ \p \llb A_i \rrb$ are oriented boundaries, each of them minimizes the perimeter in the BV-sense. There are localized versions of Proposition \ref{dic} for currents in a manifold $M^{n+1}$. When $U$ is a proper ball in $M$ and $R \in \mathcal{R}_n(U)$, we take a diffeomorphism $f:U \ra \R^{n+1}$,
apply Proposition \ref{dic} to $f_\sharp R \in \mathcal{R}_n(\R^{n+1})$ and consider the pull-back of the resulting decomposition on $U$. We state this local decomposition for area minimizers as follows.

\begin{proposition}[Local decompositions]\label{dicc} \,
Let $U\subset M$ be open with $H_n(M, M\setminus U)$, and let $T \in \mathcal{R}_n(U)$ be a locally mass minimizing current with $\p T =0$. Then there exist oriented boundaries $\p A_i$ for open $A_i\subset U$, $i \in \Z$, which are locally area minimizing in $U$ and satisfy $A_i \subset A_{i+1}$ such that for any open $W\supset U$, we have
\begin{equation}\label{ld}
T\llcorner U= \sum_{i \in \Z} \p \llb A_i \rrb \llcorner U\quad\mm{and}\quad {\bf{M}}_W (T\llcorner U)= \sum_{i \in \Z} {\bf{M}}_W (\p \llb A_i \rrb \llcorner U).
\end{equation}
\end{proposition}

The strict maximum principle~\cite[Chapter 2]{Si2} shows that the oriented boundaries in the sum~\eqref{ld} are either locally disjoint or equal for currents with multiplicities.\\

 For a mass minimizing current $T$ representing a given homology class $\alpha \in H_{n}(M^{n+1}, \Z)$ of a \emph{compact} manifold $M$, this sum is \emph{finite}. Indeed, take a small ball $B_{5r}(p) \subset U \subset M$, so that $(5r)^{-1} \cdot B_{5r}(p)$ is nearly isometric to the unit ball in $\R^{n+1}$. Then the minimality of the $\p A_i$ which intersect $B_r(p)$ gives the estimate
${\bf{M}}_{B_{2\rho}(p)}(\p \llb A_i)\rrb \llcorner U \ge c_n \cdot r^n$ for some constant $c_n >0$ only depending on $n$~\cite[Inequality (5.16)]{Gi}. Then the finiteness of the total mass of $T$ shows that there are only finitely many such $\p A_i$.\\

For this decomposition, the term \textbf{local} refers to the choice of a suitable set $U$ in the ambient manifold, \emph{independent} of the given current. This allows us to use these results, within a fixed set $U$, when we consider converging sequences of such currents.
%
%
%
%%%%%%%%%%%%%%%%%%%%%%%%%%%%%%%%%
%%%%%%%%%%%%%%%%%%%%%%%%%%%%%%%%%
%%%%%%%%%%%%%%%%%%%%%%%%%%%%%%%%%
\setcounter{section}{2}
\renewcommand{\thesubsection}{\thesection}
\subsection{\si-Whitney smoothings}\label{swi}
%%%%%%%%%%%%%%%%%%%%%%%%%%%%%%%%%
%%%%%%%%%%%%%%%%%%%%%%%%%%%%%%%%%
%%%%%%%%%%%%%%%%%%%%%%%%%%%%%%%%%
Here we explain how to define for any $H \in {\cal{H}}$ a certain locally finite ball cover of $H \setminus \Sigma$ which can be used to controllably smooth out the merely Lipschitz regular function $\delta_\bp$. The overall strategy resembles the classical Whitney smoothing in~\cite{Wh}, whence the name of \si-Whitney smoothing.\\

\noindent\textbf{I.\ Locally finite covers.} \,
We first prove the existence of Besicovich style covers of $H \setminus \Sigma$ particularly adapted to $\bp$. The proof deals with the non-totally geodesic case. For a consistent statement we also include a statement for totally geodesic $H$. There, the result boils down to the surjectivity of the exponential map and an infinitely sheeted covering if $H$ is compact.

\begin{proposition}[\si-adapted covers]\label{besi} \,
For a given \si-transform $\bp$ there exists $\xi_0\equiv\xi_0(n,L_{\bp})\in (0,1/10^3L_{\bp})$ such that for any $H \in {\cal{H}}$ and $\xi \in (0,\xi_0)$ we can construct an \textbf{\si-adapted cover} ${\cal A}$. This is a locally finite cover ${\cal{A}}=\{\overline{B}_{\Theta(p)}\,|\, p \in Z \}$ of $H\setminus\Sigma$ by closed balls of radius $\Theta(p):=\xi / \bp(p) = \xi \cdot \delta_{\bp}(p)$ with centers in a discrete set $Z\subset H \setminus \Sigma$ and such that the following properties hold: If $Q_H=Q_0\cap H \setminus \Sigma$ with $Q_0$ a suitably small neighborhood of $\Sigma$ in $H$, then
\begin{description}
  \item[(C1)] for $p \in Q_H$ the exponential map $\exp_p|_{B_{100\Theta(p)}(0)}$ is bi-Lipschitz onto its image for some bi-Lipschitz constant $l(n)\ge 1$.
  \item[(C2)] $Z^Q := Z\cap Q_H$ splits into $c(n)$ disjoint families $Z^Q(1),\ldots,Z^Q(c)$ with
  \begin{enumerate}
    \item $B_{10\Theta(p)}(p) \cap B_{10\Theta(q)}(q) \v$ for $p$ and $q$ in the same $Z^Q(k)$.
    \item $q \notin \overline{B}_{\Theta(p)}(p)$ for any two $p$, $q \in Z^Q$.
  \end{enumerate}
\end{description}
In particular, for $z \in Q_H$ and $\rho \in (0,10)$ there is a uniform bound on the covering number $\cs(Z^Q,z,\rho):= \cs\{x \in Z^Q \,|\, z \in B_{\rho \cdot \Theta(x)}(x)\} \le c(n)$. Furthermore, we have:
\begin{itemize}
  \item For any $\ve  >0$ we can find some $\xi_{\ve} \in (0,\xi_0)$ such that for every $p \in H \setminus \Sigma$ the exponential map $\exp_p|_{B_{100\xi_{\ve}/ \bp(p)}(0)}$ is  $(1 + \ve)$-bi-Lipschitz onto its image.
  \item For $H\in{\cal{H}}^\R_n$ we may choose $Q_H = H \setminus \Sigma$.
\end{itemize}
\end{proposition}

In the proof of \ref{besi} we use the following Harnack style property of $\bp$.

\begin{lemma}\label{liph} \,
Let $\bp=\bp_H>0$ and $L\equiv L(\bp)$ be the Lipschitz constant of $\delta_{\bp}$. Then $\bp$ is locally Lipschitz and
\[
\left|\bp(q)/\bp(p)  - 1\right|  \le 2L \cdot \bp(p) \cdot d_{g_H}(q,p),\mm{ for any }q \in B_{1/(2L \cdot \bp(p))}(p).
\]
\end{lemma}

\noindent{\bf Proof} \,
As reciprocal of a Lipschitz function, $\bp$ is at least locally Lipschitz. For any two $p$, $q \in H \setminus \Sigma$ the inequality $|\delta_{\bp}(p) - \delta_{\bp}(q)|\le L \cdot d_{g_H}(p,q)$ gives
\[
|\bp(p)-\bp(q)|  \le L \cdot \bp(p) \cdot \bp(q) \cdot d_{g_H}(p,q),
\]
whence $\bp(x) \le 2\bp(p)$ for all $x \in  B_{1/(2L \cdot \bp(p))}(p)$. Thus for any $q \in B_{1/(2 \cdot L \cdot \bp(p))}(p)$ we directly get $|\bp(p)-\bp(q)| \le 2 \cdot L \cdot  \bp^2(p)  \cdot d_{g_H}(p,q)$.
\qed

\noindent{\bf Proof of \ref{besi}} \,
We subdivide the proof into three steps. In the first two steps we derive pointwise estimates for the volume of balls within $B_{\Theta(p)}(p)$ using scalings. By Lemma~\ref{liph} these estimates are locally uniform. Finally we use some simple combinatorics to ensure the claimed properties of ${\cal{A}}$. Let again $L\equiv L(\bp)$ be the Lipschitz constant of $\delta_{\bp}$. %In each step we temporarily choose some local scaling of $H$ to simplify the argument but we undo these scalings before we switch to the next step.
\\

\noindent\textbf{Step 1 (Scaling of $\mathbf{\Sigma \subset H \subset M}$)}\\
For $p \in H \setminus \Sigma$ we scale $M$ by $L \cdot \bp(p)$. In particular, this produces out of $B_{1 /(L \cdot \bp(p))}(p)$ the ball $B_1(p) \subset L \cdot \bp(p) \cdot H \subset L \cdot \bp(p) \cdot M$. Since $L \cdot \bp(x) \ge 1/dist_{g_H}(x, \Sigma)$, the rescaled manifold $L \cdot \bp(p) \cdot M$ becomes virtually flat as $p$ approaches $\Sigma$. More formally, let us denote again by $\exp_p[s \cdot M]: (T_p M,g_{T_pM})\ra s \cdot M$ the exponential map of $s \cdot M$ in $p$, $s \ge 1$. Then, for any $\ve >0$, we find a neighborhood $W(\ve) \subset H$ of $\Sigma$ such that
\[
\big|\exp_p[\bp(p) \cdot M]^*(L^2 \cdot \bp(p)^2 \cdot g_M) - g_{T_p M}\big|_{C^5(B_{100}(0))} \le \ve
\]
if $p \in W(\ve)$. Here, both the $C^5$-norm and the radius are measured with respect to $\bp(p) \cdot M$.\\

\noindent\textbf{Step 2 (Locally uniform estimates on $\mathbf{H \setminus \Sigma}$)}\label{co}\\
\noindent\textbf{2.1.}\ Since we have a lower bound for $\bp|_{W(\ve)}$ which diverges as $\ve \ra 0$, rescaling of $M$ by $\bp(p)$, for $p \in W(\ve)$, shows that $M$ converges to a flat space near $p$. But approaching $\Sigma$ as $\ve\ra0$ also means that $|A|$ diverges. This time $\bp \ge |A|$ shows that $|A|(p) \le 1$ after scaling by $\bp(p)$. Lemma \ref{liph} boosts this pointwise estimate to local estimates, namely
\begin{equation}\label{aest}
\left|\bp(x)/\bp(p)  - 1\right| \le  2L \cdot \bp(p) \cdot d_{g_H}(x,p) \le 1
\end{equation}
for any $p \in H \setminus \Sigma$ and $x \in B_{1/(2L \cdot  \bp(p))}(p) \subset H$. Then \eqref{aest} and $\bp \ge |A|$ again imply that there is a constant $A_n \ge 1$ depending only on $n$ such that $|A| \le A_n$ on $B_1(p) \subset L \cdot \bp(p) \cdot H$.\\

\noindent\textbf{2.2.}\ This locally uniform bound on $|A|$ implies a locally uniform volume estimate for balls in $H$: Gauss equations relating the curvature tensors of $H$ and $M$ as well as Rauch's comparison theorem show that for some $\ve >0$ small enough, we can find positive functions $\Lambda(\zeta),\,\eta(\zeta)>0$ with
\[
\Lambda(\zeta) \ra \infty\quad\mm{and}\quad\eta(\zeta) \ra 0\mm{ as } \zeta \ra 0.
\]
For any $p \in W(\ve)$ the map $exp_p[\Lambda(\zeta) \cdot  L \cdot \bp(p)\cdot H ]$ is therefore a local diffeomorphism from $B_{10^3}(0)$ onto its image in $H$ with
\[
\big|\exp_p[\Lambda(\zeta) \cdot  L \cdot \bp(p)\cdot H ]^*(\Lambda^2(\zeta) \cdot  L^2 \cdot \bp(p)^2 \cdot g_H) - g_{T_p H}\big|_{L^\infty(B_{100}(0))} \le \eta(\zeta)
\]
measured with respect to $\Lambda(\zeta) \cdot L \cdot \bp(p) \cdot H$. By the regularity theory of $H$ this can be upgraded to $C^k$-estimates for any given $k \ge 0$. We then obtain, keeping the same notation for the constants for simplicity,
\[
\big|\exp_p[\Lambda(\zeta) \cdot  L \cdot \bp(p) \cdot H ]^*(\Lambda^2(\zeta)    \cdot L^2 \cdot \bp(p)^2 \cdot g_H) - g_{T_p H}\big|_{C^5(B_{100}(0))} \le\eta(\zeta).
\]
\noindent\textbf{2.3.}\ Choosing $\ve >0$ sufficiently small we acquire uniform control for any $p \in W(\ve)$:\\
\noindent2.3A.\ The map $\exp_p[\Lambda(\zeta) \cdot  L \cdot \bp(p) \cdot H]$ is bi-Lipschitz from $B_{100}(p)$ to its image for some bi-Lipschitz constant $l(\zeta) \ge 1$ with $l(\zeta) \ra 1$ for $\zeta \ra 0$.\\
\noindent2.3B.\ If $\zeta >0$ so that $l(\zeta) \in [1,2]$, then for $z \in B_{50}(p)$ the volume estimates
\begin{equation}\label{vol}
k_1 \le Vol(B_{1/3}(z))\mm{ and }Vol(B_{30}(z))  \le k_2, \mm{ for constants }k_i(n,\bp) >0,\,i=1,\,2
\end{equation}
hold (with volumes and radii taken with respect to $\Lambda(\zeta) \cdot   L \cdot \bp(p) \cdot H$).\\

\noindent\textbf{Step 3 (Combinatorics)}\\
Pick some $\zeta \ll 1$ so that $\Lambda(\zeta) \gg 1$ and $L \cdot \Lambda(\zeta) >100$. We set
\[
\xi(\zeta):=1/(L \cdot \Lambda(\zeta))\quad\mm{and}\quad\Theta(p):= \xi(\zeta) / \bp(p) = 1/(L \cdot \Lambda(\zeta) \cdot \bp(p)).
\]
Next we choose a countable dense subset $S$ of $H \setminus \Sigma$, $S= \{a_m\,|\,m \in \Z^{\ge 0}\}$ and define a cover by ${\cal{B}} = \{\overline{B}_{\Theta(p)(p)}\,|\, p \in S\}$. We define a map $i: S \ra\Z^{\ge 0}$ by induction: Set $i(a_0):=1$ and
\[
i(a_{k+1}):=\left\{\begin{array}{l}0, \mm{ if } a_{k+1} \in \bigcup_{i \le k} \overline{B}_{\Theta(a_i)(a_i)}\\[5pt]
\min\big(\{m \le k \,|\,d_{g_H}(a_m,a_{k+1})/10 > \Theta(a_m) + \Theta(a_{k+1})\} \cup \{k+1\}\big), \mm{ otherwise.}\end{array}\right.
\]
Finally, we put
\[
{\cal{A}} := \bigcup_{j\ge 1} {\cal{A}}(j) \mm{ with }{\cal{A}}(j) := \{\overline{B}_{\Theta(p)}(p)\,|\, p \in S,\,i(p)=j\},\, \mm{ and } Z:= \{a \in S\,|\, i(a) \ge 1\}.
\]
Then these families satisfy ${\cal{A}}(i) \cap {\cal{A}}(j) \v$ for $i \neq j$ and also (i) and (ii) of (C2). Moreover, there is a neighborhood $Q_0$ of $\Sigma$ and a constant $c\equiv c(n,\bp)$ such that ${\cal{A}}(i) \v$ for $i > c$. In particular, this implies the local finiteness of ${\cal{A}}$ since
there are only finitely many balls in ${\cal{A}}$ with center in $H \setminus Q_0$. From \eqref{aest} we may assume that
\[
9/10 \cdot \Theta(p) \le  \Theta(x) \le 11/10 \cdot \Theta(p)\mm{ for any }x \in \overline{B}_{100 \cdot \Theta(p)}(p).
\]
We put $c(n,\bp):= \mm{\emph{the smallest integer which is}}\ge k_2/k_1$ with $k_i$ as in \eqref{vol}, and claim that ${\cal{A}}(i) \v$ for $i > c$. Otherwise, we could take $\overline{B}_{\Theta(p)}(p)\in {\cal{A}}(c+1)$, whence $B_{10 \cdot \Theta(p)}(p) \cap B_{10 \cdot \Theta(x_i)}(x_i) \n$ for at least $c$ different $x_i \in Z$. Now we use \eqref{vol} with respect to $\Lambda(\zeta) \cdot   L \cdot \bp(p) \cdot H$ and obtain
\begin{align*}
(c+1) \cdot k_1 &\le \sum Vol(B_{1/3}(x_i))+ Vol(B_{1/3}(p))\\
&= Vol \big(\bigcup B_{1/3}(x_i) \cup B_{1/3}(p)\big) \le Vol(B_{30}(p)) \le k_2.
\end{align*}
But then the $(c+1)$ balls $B_{1/3}(x_i)$ and $B_{1/3}(p)$ are pairwise disjoint, a contradiction. Finally, we observe that ${\cal{A}}$ is indeed a cover. Assume there is some point $q \in U=H\setminus\Sigma \setminus \bigcup_{p \in A} \overline{B}_{\Theta(p)}(p)$. Then $U$ is open since ${\cal{A}}$ is locally finite, so that there is a point $z \in U\cap S$ with $\overline{B}_{\Theta(z)}(z) \in {\cal{A}}$ and $q \in \overline{B}_{\Theta(z)}(z)$, a contradiction.
\qed

\noindent\textbf{II.\ Whitney smoothings.} \,
In $\R^n$, the metric distance to a closed subset is a Lipschitz function. Whitney introduced in~\cite{Wh} a method to smooth out the distance function while keeping most of the information it carries. We mimick his proof using \si-adapted covers.

\begin{proposition}[\si-Whitney smoothings]\label{smsk} \,
For any \si-transform $\bp$ there is smoothing $\bp^*$, i.e., a family of smooth functions $\bp_H^*$ defined on $H \setminus \Sigma$ for any $H \in {\cal{H}}$. $\bp^*$ still satisfies axioms (S1) - (S3) for \si-transforms and we have:
\begin{equation}\label{smot}
c_1 \cdot \delta_{\bp}(x) \le \delta_{\bp^*}(x)  \le c_2 \cdot \delta_{\bp}(x)\quad\mm{and}\quad |\p^\beta \delta_{\bp^*}  / \p x^\beta |(x) \le c_3(\beta) \cdot \delta_{\bp}^{1-|\beta|}(x)
\end{equation}
for constants $c_i >0$, $i=1,\,2,\,3$. Here, $\beta$ is a multi-index for derivatives with respect to normal coordinates around $x \in H \setminus \Sigma$. We have the following dependancies.
\begin{enumerate}
  \item For $H \in {\cal{H}}_n$, we have on $Q_H \subset H\setminus\Sigma$: $c_{1,2}=c_{1,2}(L_{\bp}, n)$, and $c_3=c_3(L_{\bp},n,\beta)$.
  \item For $H \in {\cal{H}}^{\R}_n$, we have on $H \setminus \Sigma$: $c_{1,2}=c_{1,2}(L_{\bp},n)$, and $c_3=c_3(L_{\bp},n,\beta)$.
  \item For $H \in {\cal{H}}^c_n$, we have on $H \setminus \Sigma$: $c_{1,2}=c_{1,2}(L_{\bp},H)$, and $c_3=c_3(L_{\bp},H,\beta)$.
\end{enumerate}
We interpret this as some weakened naturality of $\bp^*$, we call its quasi-naturality.
\end{proposition}

\noindent{\bf Proof } \,
For $\bp_H \equiv 0$, we set $\bp^*=0$; otherwise, we have $\bp>0$. We choose a smooth non-negative function $\phi$ on $\R^n$ with $\phi \equiv 1$ on $B_1(0)$ and $\phi \equiv 0$ on $\R^n  \setminus B_2(0)$. As in Proposition~\ref{besi} Step 2.2.\ we consider $\exp_p[10 \cdot \Lambda(\zeta) \cdot \bp(p) \cdot H]$, for some sufficiently small $\zeta >0$, and define $\Phi_{p}(x):= \phi\big(\exp^{-1}_p[10\Lambda(\zeta) \cdot \bp(p) \cdot H](x)\big)$, $x \in H \setminus \Sigma$. Thus the $exp_p$-preimage of the ball $B_5(p) \subset 10\Lambda(\zeta) \cdot \bp(p) \cdot H$ is almost isometric to $B_5(0) \subset T_pH$. We notice that $|\p^\beta \Phi_p/ \p x^\beta |(x)\le k(\beta) \cdot (10\Lambda(\zeta) \cdot \bp(x))^{|\beta|-1}$ on $B_5(p) \subset 10\Lambda(\zeta) \cdot \bp(p)
\cdot H$. Now we define $\bp^*$ through its \si-distance, namely
\[
\delta_{\bp^*}(x):= \sum_{p \in A} \delta_{\bp}(p) \cdot \Phi_{p}(x).
\]
This is a smooth positive function. To check \eqref{smot} and (i)-(iii) we note that (i) implies (ii) and (iii) by Proposition~\ref{besi}. On the other hand, (i) readily follows from Proposition~\ref{besi} for sufficiently small $\ve>0$ by using the upper bound on the  covering number of ${\cal{A}}$, the $(1+\ve)$-Lipschitz charts of Proposition~\ref{besi} (ii), and the Lipschitz condition on $\delta_{\bp}$.
\qed

\end{document}